\documentclass[11pt]{article}
\usepackage{epic,latexsym,amssymb}
\usepackage{color}
\usepackage{tikz}
\usepackage{amsfonts,epsf,amsmath,leftidx}

\newtheorem{Theorem}{Theorem}
\newtheorem{Lemma}{Lemma}

\newtheorem{Corollary}{Corollary}
\newtheorem{ob}{Observation}

\newcommand{\qed}{$\Box$}

\newcommand{\cP}{\mathcal{P}}

\newcommand{\cC}{{\cal C}}

\newcommand{\proof}{\noindent\textbf{Proof. }}

\let\oldenumerate\enumerate
\renewcommand{\enumerate}{
  \oldenumerate
  \setlength{\itemsep}{0pt}
  \setlength{\parskip}{0pt}
  \setlength{\parsep}{0pt}
}

\begin{document}

\title{Structural Properties of Connected Domination Critical Graphs}


\author{Norah Almalki$^1$ \, 
and  Pawaton Kaemawichanurat$^{2,3}$\thanks{Research
supported by Thailand Research Fund (MRG 6080139 and MRG 6280223)}
\\ \\
$^{1}$Department of Mathematics and Statistics \\
Taif University, Taif City \\
Saudi Arabia\\
\small \tt Email: norah@tu.edu.sa \\
\\
$^{2}$Department of Mathematics, Faculty of Science \\
King Mongkut's University of Technology Thonburi\\
Bangkok, Thailand \\
\small \tt Email: pawaton.kae@kmutt.ac.th \\
\\
$^{3}$Mathematics and Statistics with Applications~(MaSA)}

\date{}
\maketitle


\begin{abstract}
A graph $G$ is said to be $k$-$\gamma_{c}$-critical if the connected domination number $\gamma_{c}(G)$ is equal to $k$ and $\gamma_{c}(G + uv) < k$ for any pair of non-adjacent vertices $u$ and $v$ of $G$. Let $\zeta$ be the number of cut vertices of $G$ and let $\zeta_{0}$ be the maximum number of cut vertices that can be contained in one block. For an integer $\ell \geq 0$, a graph $G$ is $\ell$-factor critical if $G - S$ has a perfect matching for any subset $S$ of vertices of size $\ell$. It was proved that, for $k \geq 3$, every $k$-$\gamma_{c}$-critical graph has at most $k - 2$ cut vertices and the graphs with maximum number of cut vertices were characterized. It was proved further that, for $k \geq 4$, every $k$-$\gamma_{c}$-critical graphs satisfies the inequality $\zeta_{0}(G) \le \min \left\{ \left\lfloor \frac{k + 2}{3} \right\rfloor, \zeta \right\}$. In this paper, we characterize all $k$-$\gamma_{c}$-critical graphs having $k - 3$ cut vertices. Further, we establish realizability that, for given $k \geq 4$, $2 \leq \zeta \leq k - 2$ and $2 \leq \zeta_{0} \le \min \left\{ \left\lfloor \frac{k + 2}{3} \right\rfloor, \zeta \right\}$, there exists a $k$-$\gamma_{c}$-critical graph with $\zeta$ cut vertices having a block which contains $\zeta_{0}$ cut vertices. Finally, we proved that every $k$-$\gamma_{c}$-critical graph of odd order with minimum degree two is $1$-factor critical if and only if $1 \leq k \leq 2$. Further, we proved that every $k$-$\gamma_{c}$-critical $K_{1, 3}$-free graph of even order with minimum degree three is $2$-factor critical if and only if $1 \leq k \leq 2$.
\end{abstract}
{\small \textbf{Keywords:}Domination; Characterization; Matching; Realizability


\section{\bf Introduction}
For a natural number $n$, we let $[n] = \{1, 2, ..., n\}$. All graphs in this paper are finite, undirected and simple (no loops or multiple edges). For a graph $G$, let $V(G)$ denote the set of vertices of $G$ and let $E(G)$ denote the set of edges of $G$. For $S \subseteq V(G)$, $G[S]$ denotes the subgraph of $G$ induced by $S$. The \emph{open neighborhood} $N_{G}(v)$ of a vertex $v$ in $G$ is $\{u \in V(G) : uv \in E(G)\}$. Further, the \emph{closed neighborhood} $N_{G}[v]$ of a vertex $v$ in $G$ is $N_{G}(v) \cup \{v\}$. For subsets $X$ and $Y$ of $V(G)$, $N_{Y}(X)$ is the set $\{y \in Y : yx \in E(G)$ for some $x \in X \}$. For a subgraph $H$ of $G$, we use $N_{Y}(H)$ instead of $N_{Y}(V(H))$ and we use $N_{H}(X)$ instead of $N_{V(H)}(X)$. If $X = \{x\}$, we use $N_{Y}(x)$ instead of $N_{Y}(\{x\})$. The \emph{degree} $deg(x)$ of a vertex $x$ in $G$ is $|N_{G}(x)|$. The \emph{minimum degree} of $G$ is denoted by $\delta(G)$. When no ambiguity occur, we write $N(x), N(X)$ and $\delta(G)$ instead of $N_{G}(x), N_{G}(X)$ and $\delta$, respectively. An \emph{end vertex} is a vertex of degree one and a \emph{support vertex} is the vertex which is adjacent to an end vertex. A \emph{tree} is a connected graph with no cycle. A \emph{star} $K_{1, n}$ is a tree containing one support vertex and $n$ end vertices. The support vertex of a star is called the \emph{center}. For a connected graph $G$, a vertex subset $S$ of $G$ is called a \emph{cut set} if $G - S$ is not connected. We let $\omega_{o}(G - S)$ be the number of components of $G - S$ containing odd number of vertices. In particular, if $S = \{v\}$, then $v$ is called a \emph{cut vertex} of $G$. That is, $G - v$ is not connected. A \emph{block} $B$ of a graph $G$ is a maximal connected subgraph such that $B$ has no cut vertex. An \emph{end block} of $G$ is a block containing exactly one cut vertex of $G$. For graphs $H$ and $G$, a graph $G$ is said to be $H$-\emph{free} if $G$ does not contain $H$ as an induced subgraph. For a connected graph $G$, a \emph{bridge} $xy$ of $G$ is an edge such that $G - xy$ is not connected.
\vskip 5 pt

\indent For a finite sequence of graphs $G_{1}, ..., G_{l}$ for $l \geq 2$, the \emph{joins} $G_{1} \vee  \cdots \vee G_{l}$ is the graph consisting of the disjoint union of $G_{1}, ..., G_{l}$ and joining edges from every vertex of $G_{i}$ to every vertex of $G_{i + 1}$ for $1 \leq i \leq l - 1$. In particular, for a subgraph $H$ of $G_{2}$, the \emph{join} $G_{1} \vee \leftidx{_H}{G}{_{2}}$ is the graph consisting of the disjoint union of $G_{1}$ and $G_{2}$ and joining edges from every vertex of $G_{1}$ to every vertex of $H$. As the join operation is run over vertices, for a vertex $x$ and a set $X$ of vertices, the \emph{join} $x \vee X$ is the graph consisting of the disjoint union of $\{x\}$ and $X$ and joining edges from $x$ to every vertex in $X$. For a subgraph $H$ of $G$, $H$ is \emph{maximal complete subgraph} of $G$ if for any complete subgraph $H'$ of $G$ such that $V(H) \subseteq V(H')$, we have $V(H) = V(H')$.
\vskip 5 pt

\indent The \emph{distance} $d(u, v)$ between vertices $u$ and $v$ of $G$ is the length of a shortest $(u, v)$-path in $G$. The \emph{diameter} of $G$ $diam(G)$ is the maximum distance of any two vertices of $G$. For a non-negative integer $k$, a graph $G$ is $k$-\emph{diameter critical} if $diam(G) = k$ and $diam(G - uv) > k$ for any edge $uv \in E(G)$. A \emph{matching} of a graph $G$ is a set of edges which are not incident to a common vertex. A matching $M$ of a graph $G$ is \emph{perfect} if $V(M) = V(G)$. For a non-negative integer $\ell$, a graph $G$ is $\ell$-\emph{factor critical} if $G - S$ has a perfect matching where $S$ is any set of $\ell$ vertices, in particular, a graph $G$ is \emph{factor critical} if $\ell = 1$ and is \emph{bi-critical} if $\ell = 2$. For subsets $D$ and $X$ of $V(G)$, $D$ \emph{dominates} $X$ if every vertex in $X$ is either in $D$ or adjacent to a vertex in $D$. If $D$ dominates $X$, then we write $D \succ X$ and we also write $a \succ X$ when $D = \{a\}$. Moreover, if $X = V(G)$, then $D$ is a \emph{dominating set} of $G$ and we write $D \succ G$ instead of $D \succ V(G)$. A \emph{connected dominating set} of a graph $G$ is a dominating set $D$ of $G$ such that $G[D]$ is connected. If $D$ is a connected dominating set of $G$, we then write $D \succ_{c} G$. A smallest connected dominating set is called a $\gamma_{c}$\emph{-set}. The cardinality of a $\gamma_{c}$-set is called the \emph{connected domination number} of $G$ and is denoted by $\gamma_{c}(G)$. A graph $G$ is said to be $k$-$\gamma_{c}$\emph{-critical} if $\gamma_{c}(G) = k$ and $\gamma_{c}(G + uv) < k$ for any pair of non-adjacent vertices $u$ and $v$ of $G$.
\vskip 5 pt

\indent In the structural characterizations of $k$-$\gamma_{c}$-critical graphs, Chen et al.\cite{XGC} showed that every $1$-$\gamma_{c}$-critical graph is a complete graph while every $2$-$\gamma_{c}$-critical graph is the complement of the disjoint union of at least two stars. However, for $k = 3$, it turns out that the $k$-$\gamma_{c}$-critical graphs have no complete characterization in the sense of free graphs (see \cite{KC}, Chapter 5). Interestingly, it was proved by Hanson and Wang \cite{HW} that, for a connected graph $G$, the graph $G$ is $3$-$\gamma_{c}$-critical if and only if the complement of $G$ is $2$-diameter critical. For a study on $k$-diameter critical graphs see Almalki \cite{NC}. In \cite{A}, Ananchuen proved that every $3$-$\gamma_{c}$-critical graph contains at most one cut vertex and also established characterizations of $3$-$\gamma_{c}$-critical graphs having a cut vertex. For more studies related with $3$-$\gamma_{c}$-critical graphs see \cite{AW,LP08,MA}. For $k = 4$, Kaemawichanurat and Ananchuen \cite{PKNA} proved that every $4$-$\gamma_{c}$-critical graph contains at most two cut vertices and the characterization of the such graphs having two cut vertices was given. Further, Kaemawichanurat and Ananchuen \cite{PKNA2} established that every $k$-$\gamma_{c}$-critical graphs contains at most $k - 2$ cut vertices when $k \geq 5$. They also characterized that there is exactly one class of $k$-$\gamma_{c}$-critical graphs having $k - 2$ cut vertices. In the same paper, the authors established the maximum number of cut vertices that every block of the graph can have. That is:

\begin{Theorem}{\rm (\cite{PKNA2})}
\label{thm mp2}
Let $G$ be a $k$-$\gamma_{c}$-critical graph containing $\zeta$ cut vertices and let $\zeta_{0}(G)$ be the maximum number of cut vertices of $G$ that can be in a block of $G$. Then
\[
\zeta_{0}(G) \le \min \left\{ \left\lfloor \frac{k + 2}{3} \right\rfloor, \zeta \right\}.
\]
\end{Theorem}

\noindent Very recently, Henning and Kaemawichanurat \cite{HK} characterized all the eleven classes of $k$-$\gamma_{c}$-critical graphs satisfying the upper bound of Theorem \ref{thm mp2}.
\vskip 5 pt

\indent In this paper, for $k \geq 5$, we characterize all $k$-$\gamma_{c}$-critical graphs having $k - 3$ cut vertices. Further, we establish realizability that, for given $k \geq 4$, $2 \leq \zeta \leq k - 2$ and $2 \leq \zeta_{0} \le \min \left\{ \left\lfloor \frac{k + 2}{3} \right\rfloor, \zeta \right\}$, there exists a $k$-$\gamma_{c}$-critical graph with $\zeta$ cut vertices having a block which contains $\zeta_{0}$ cut vertices. We also proved that every $k$-$\gamma_{c}$-critical graph of odd order with $\delta \geq 2$ is factor critical  if and only if $1 \leq k \leq 2$. We prove that every $k$-$\gamma_{c}$-critical $K_{1, 3}$-free graph of even order with $\delta \geq 2$ is factor critical  if and only if $1 \leq k \leq 2$. All the main results are shown in Section \ref{main} while their proofs are given in Section \ref{proof}. We present some useful results that are used in the proofs in Section \ref{pre}.

\vskip 5 pt

\section{Main Results}\label{main}
\noindent First, we give constructions of two classes of such graphs which will be the characterizations of $k$-$\gamma_{c}$-critical graphs with $k - 3$ cut vertices. Let
\vskip 5 pt

\indent $\textbf{\emph{i}} = (i_{1}, i_{2}, ..., i_{k - 3})$
\vskip 5 pt

\noindent be a $k - 3$ tuples such that $i_{1}, i_{2}, ..., i_{k - 3} \in \{0, 1\}$ and $\Sigma^{k - 3}_{j = 1}i_{j} = 1$ (there is exactly one $l \in \{1, 2, ..., k - 3\}$ such that $i_{l} = 1$ and $i_{l'} = 0$ for all $l' \in \{1, 2, ..., k - 3\} - \{l\}$).
\vskip 18 pt

\noindent \textbf{The class} $\mathcal{G}_{1}(i_{1}, i_{2}, ..., i_{k - 3})$\\
\indent For a $k - 3$ tuples $\textbf{\emph{i}} = (0, 0, ..., i_{l}, ..., 0)$ where $i_{l} = 1$ and $i_{l'} = 0$ for $1 \leq l \neq l' \leq k - 2$, a graph $G$ in the class $\mathcal{G}_{1}\textbf{\emph{i}}$ can be constructed from paths $c_{0}, c_{1}, ..., c_{l - 1}$ and $c_{l}, c_{l + 1}, ..., c_{k - 4}$, a copy of a complete graph $K_{n_{l}}$ and a block $B \in \mathcal{B}_{2, 2}$ (the construction of the class $\mathcal{B}_{2, 2}$ will be given in Section \ref{pre} as it was established earlier in \cite{PKNA2}) by adding edges according the join operations :
\begin{itemize}
  \item $c_{l - 1} \vee K_{n_{l}} \vee c_{l}$ and
  \item $c_{k - 4} \vee c$
\end{itemize}
\noindent where we call $c$ the \emph{head} of $B$. Examples of graphs in this case are illustrated by Figures \ref{k-3-1} and \ref{k-3-2}. Further, for a $k - 3$ tuples $\textbf{\emph{i}} = (0, 0, ..., 1)$ where $i_{k - 3} = 1$ and $i_{l'} = 0$ for $1 \leq l' \leq k - 2$, a graph $G$ in the class $\mathcal{G}_{1}\textbf{\emph{i}}$ can be constructed from paths $c_{0}, c_{1}, ..., c_{k - 4}$, a copy of a complete graph $K_{n_{k - 3}}$ and a block $B \in \mathcal{B}_{2, 2}$ by adding edges according the join operation :
\begin{itemize}
  \item $c_{k - 4} \vee K_{n_{k - 3}} \vee c$
\end{itemize}
\noindent where we call $c$ the \emph{head} of $B$. Example of a graph in this case is illustrated by Figure \ref{k-3-3}.
\vskip 5 pt

\begin{figure}[htb]
\begin{center}
\setlength{\unitlength}{0.8cm}
\begin{picture}(15, 5)
\put(3, 3){\circle*{0.2}}
\put(4.2, 3){\oval(1, 2)}
\put(11, 3){\oval(3, 2)}

\put(3, 3){\line(1, 1){0.7}}
\put(3, 3){\line(1, -1){0.7}}
\put(5.4, 3){\line(-1, 1){0.7}}
\put(5.4, 3){\line(-1, -1){0.7}}

\put(5.4, 3){\line(1, 0){1}}
\put(8, 3){\line(1, 0){1}}

\put(9, 3){\line(1, 1){0.7}}
\put(9, 3){\line(1, -1){0.7}}

\put(9, 1.5){\line(1, 0){3.5}}
\put(9, 1.5){\line(0, 1){0.5}}
\put(12.5, 1.5){\line(0, 1){0.5}}
\put(10.75, 1.5){\line(0, -1){0.5}}

\put(5.4, 3){\circle*{0.2}}
\put(6.4, 3){\circle*{0.2}}
\put(8, 3){\circle*{0.2}}
\put(9, 3){\circle*{0.2}}
\put(7, 3){\footnotesize$...$}
\put(2.8, 2.3){\footnotesize$c_{0}$}
\put(5.3, 2.3){\footnotesize$c_{1}$}
\put(6.3, 2.3){\footnotesize$c_{2}$}
\put(7.7, 2.3){\footnotesize$c_{k - 4}$}
\put(8.9, 2.3){\footnotesize$c$}

\put(4, 1.5){\footnotesize$K_{n_{1}}$}
\put(10, 0.5){\footnotesize$B \in \mathcal{B}_{2, 2}$}

\end{picture}
\caption{A graph $G$ in the class $\mathcal{G}_{1}(1, 0, 0, ..., 0)$}
\label{k-3-1}
\end{center}
\end{figure}
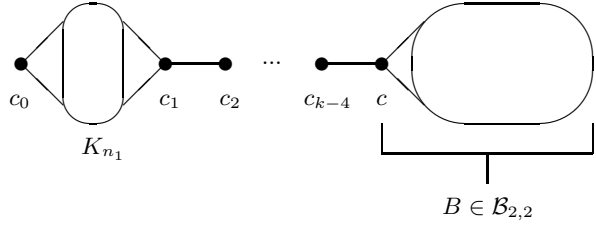
\vskip 15 pt

\begin{figure}[htb]
\begin{center}
\setlength{\unitlength}{0.8cm}
\begin{picture}(12, 4)
\put(2.5, 3){\circle*{0.2}}
\put(1.5, 3){\circle*{0.2}}
\put(0.5, 3){\circle*{0.2}}

\put(5.2, 3){\oval(1, 2)}
\put(11, 3){\oval(3, 2)}

\put(0.5, 3){\line(1, 0){2}}
\put(4, 3){\line(1, 1){0.7}}
\put(4, 3){\line(1, -1){0.7}}
\put(6.4, 3){\line(-1, 1){0.7}}
\put(6.4, 3){\line(-1, -1){0.7}}

\put(8, 3){\line(1, 0){1}}
\put(9, 3){\line(1, 1){0.7}}
\put(9, 3){\line(1, -1){0.7}}

\put(9, 1.5){\line(1, 0){3.5}}
\put(9, 1.5){\line(0, 1){0.5}}
\put(12.5, 1.5){\line(0, 1){0.5}}
\put(10.75, 1.5){\line(0, -1){0.5}}

\put(4, 3){\circle*{0.2}}
\put(6.4, 3){\circle*{0.2}}
\put(8, 3){\circle*{0.2}}
\put(9, 3){\circle*{0.2}}
\put(7, 3){\footnotesize$...$}
\put(3.1, 3){\footnotesize$...$}

\put(0.5, 2.3){\footnotesize$c_{0}$}
\put(1.5, 2.3){\footnotesize$c_{1}$}

\put(3.8, 2.3){\footnotesize$c_{l - 1}$}
\put(6.3, 2.3){\footnotesize$c_{l}$}
\put(7.7, 2.3){\footnotesize$c_{k - 4}$}
\put(8.9, 2.3){\footnotesize$c$}

\put(5, 1.5){\footnotesize$K_{n_{l}}$}
\put(10, 0.5){\footnotesize$B \in \mathcal{B}_{2, 2}$}

\end{picture}
\caption{A graph $G$ in the class $\mathcal{G}_{1}(0, 0, ..., i_{l} = 1, 0, ..., 0)$}
\label{k-3-2}
\end{center}
\end{figure}
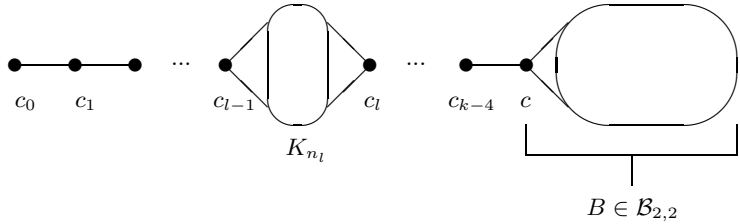
\vskip 15 pt

\begin{figure}[htb]
\begin{center}
\setlength{\unitlength}{0.8cm}
\begin{picture}(11, 4)
\put(2.5, 3){\circle*{0.2}}
\put(1.5, 3){\circle*{0.2}}
\put(0.5, 3){\circle*{0.2}}

\put(5.2, 3){\oval(1, 2)}
\put(8.3, 3){\oval(3, 2)}

\put(0.5, 3){\line(1, 0){2}}
\put(4, 3){\line(1, 1){0.7}}
\put(4, 3){\line(1, -1){0.7}}
\put(6.4, 3){\line(-1, 1){0.7}}
\put(6.4, 3){\line(-1, -1){0.7}}

\put(6.3, 3){\line(1, 1){0.7}}
\put(6.3, 3){\line(1, -1){0.7}}

\put(6.3, 1.5){\line(1, 0){3.5}}
\put(6.3, 1.5){\line(0, 1){0.5}}
\put(9.8, 1.5){\line(0, 1){0.5}}
\put(8.05, 1.5){\line(0, -1){0.5}}

\put(4, 3){\circle*{0.2}}
\put(6.4, 3){\circle*{0.2}}
\put(3.1, 3){\footnotesize$...$}

\put(0.5, 2.3){\footnotesize$c_{0}$}
\put(1.5, 2.3){\footnotesize$c_{1}$}

\put(3.8, 2.3){\footnotesize$c_{k - 4}$}
\put(6.3, 2.3){\footnotesize$c$}

\put(5, 1.5){\footnotesize$K_{n_{k - 3}}$}
\put(7.3, 0.5){\footnotesize$B \in \mathcal{B}_{2, 2}$}

\end{picture}
\caption{A graph $G$ in the class $\mathcal{G}_{1}(0, 0, ..., 1)$}
\label{k-3-3}
\end{center}
\end{figure}
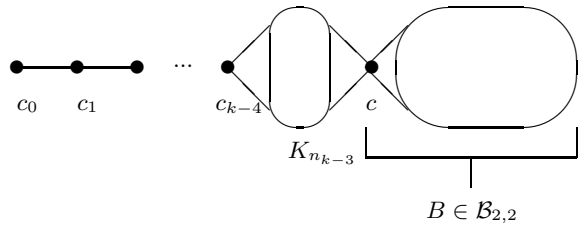

\vskip 20 pt

\indent Next, we will construct another class of $k$-$\gamma_{c}$-critical graphs with $k - 3$ cut vertices. Before giving the construction, we introduce the following class of end blocks.
\vskip 5 pt

\noindent \textbf{The class} $\mathcal{B}_{3}$\\
\indent An end block $B \in \mathcal{B}_{3}$ has $b$ as the head. Let $N_{B}(b) = A$ and $\check{B} = G[V(B) - \{b\}]$. Moreover, $B$ has the following properties
\begin{enumerate}
  \item [(1)] Every vertex $v \in V(\check{B})$, there exists a $\gamma_{c}$-set $D_{v}$ of $B$ of size $3$ such that $v \in D_{v}$.
  \item [(2)] For every non-adjacent vertices $x$ and $y$ of $\check{B}$, there exists a $\gamma_{c}$-set $D^{B}_{xy}$ of $B + xy$ such that $D^{B}_{xy} \cap \{x, y\} \neq \emptyset$, $|D^{B}_{xy}| = 2$ and $D^{B}_{xy} \cap A \neq \emptyset$.
\end{enumerate}
\vskip 5 pt

\indent It is worth noting that $D_{v}$ in the property $(1)$ satisfies $D_{v} \cap A \neq \emptyset$ in order to dominate $b$. We are ready to give the construction.
\vskip 5 pt

\noindent \textbf{The class} $\mathcal{G}_{2}(k)$\\
\indent  For $k \geq 5$, a graph $G$ in this class can be constructed from a path $c_{0}, c_{1}, ..., c_{k - 4}$ and an end block $B \in \mathcal{B}_{3}$ with the head $b$ by adding the edge $c_{k - 4}b$. For the sake of convenience, we may relabel $b$ as $c_{k - 3}$.
\vskip 8 pt

\noindent Then, we let
\vskip 8 pt

\indent $\mathcal{Z}(k, \zeta) :$ the class of $k$-$\gamma_{c}$-critical graphs containing $\zeta$ cut vertices.
\vskip 8 pt

\noindent Our first main result is:
\vskip 5 pt

\begin{Theorem}\label{thm k-3}
For an integer $k \geq 4$, $\mathcal{Z}(k, k - 3) = \mathcal{G}_{1}(i_{1}, i_{2}, ..., i_{k - 3}) \cup \mathcal{G}_{2}(k)$.
\end{Theorem}
\vskip 5 pt

\noindent For our next main results, we let $B$ be a block of $G$. We further define the following notation. We let $\cC(G)$ be the set of all cut vertices of $G$. That is
\vskip 5 pt

\indent $\cC(G) = \{ v \in V(G) : v$ is a cut vertex of $G\}$,
\vskip 5 pt

\indent $\zeta(G) = |\cC(G)|$,
\vskip 5 pt

\indent $\cC(B) = V(B) \cap \cC(G)$,
\vskip 5 pt

\indent $\zeta(B) = |\cC(B)|$ and
\vskip 5 pt

$\zeta_{0}(G) = \max\{\zeta(B) : B$ is a block of $G\}$.
\vskip 5 pt

\noindent We prove that :

\begin{Theorem}
\label{thm:main111}
For all $k \geq 4$, $2 \leq \zeta \leq k - 2$ and $2 \leq \zeta_{0} \le \min \left\{ \left\lfloor \frac{k + 2}{3} \right\rfloor, \zeta \right\}$, there exists a $k$-$\gamma_{c}$-critical graph with $\zeta$ cut vertices having a block that contains $\zeta_{0}$ cut vertices.
\end{Theorem}
\vskip 5 pt

\noindent Finally, we establish a constructive proofs to show that:
\vskip 5 pt

\begin{Theorem}\label{thm 0}
Every $k$-$\gamma_{c}$-critical graph of odd order with $\delta \geq 2$ is factor critical if and only if $k \in [2]$.
\end{Theorem}

\begin{Theorem}\label{thm bicri}
Every $k$-$\gamma_{c}$-critical $K_{1, 3}$-free graph of even order with $\delta \geq 3$ is bi-critical if and only if $k \in [2]$.
\end{Theorem}

\vskip 12 pt

\section{Preliminaries}\label{pre}
In this section, we state a number of results that are used in establishing our theorems. We begin with a result of Favaron~\cite{F} which gives matching properties of graphs according to the toughness.
\vskip 5 pt

\begin{Theorem}~\cite{F}\label{thm favaron}
For an integer $\ell \geq 0$, let $G$ be a graph with minimum degree $\delta \geq \ell + 1$. Then $G$ is $\ell$-factor critical if and only if $\omega_{o}(G - S) \leq |S| - \ell$ for any cut set $S$ of $G$ such that $|S| \geq \ell$.
\end{Theorem}
\vskip 5 pt

\indent In the context of $k$-$\gamma_{c}$-critical graphs, Ananchuen et al.~\cite{A2} established some matching properties of the such graphs when $k = 3$.
\vskip 5 pt

\begin{Theorem}\label{thm anan}\cite{A2}
Every $3$-$\gamma_{c}$-critical graph of even order with $\delta \geq 2$ has a perfect matching.
\end{Theorem}
\vskip 5 pt

\noindent However, for factor criticality, they found a $3$-$\gamma_{c}$-critical graph of odd order with $\delta \geq 2$ which is not factor critical. The graph is constructed from a complete graph $K_{n}$, a star $K_{1, n}$ by joining every end vertex of the star to every vertex of $K_{n}$ and then remove all edges of a perfect matching between these two graphs. The resulting graph is detailed in Figure~\ref{f:fig001}.

\vskip 20 pt
\begin{figure}[htb]
\begin{center}
\setlength{\unitlength}{0.7cm}
\begin{picture}(10, 5)
\put(1, 3){\circle*{0.2}}
\put(3, 3){\circle*{0.2}}
\put(3, 4){\circle*{0.2}}
\put(3, 5){\circle*{0.2}}
\put(3, 1){\circle*{0.2}}

\put(6, 3){\circle*{0.2}}
\put(6, 4){\circle*{0.2}}
\put(6, 5){\circle*{0.2}}
\put(6, 1){\circle*{0.2}}

\put(3, 1.7){\circle*{0.08}}
\put(3, 2){\circle*{0.08}}
\put(3, 2.3){\circle*{0.08}}

\put(6, 1.7){\circle*{0.08}}
\put(6, 2){\circle*{0.08}}
\put(6, 2.3){\circle*{0.08}}

\multiput(3, 3)(0.4,0){8}{\line(1,0){0.2}}
\multiput(3, 4)(0.4,0){8}{\line(1,0){0.2}}
\multiput(3, 5)(0.4,0){8}{\line(1,0){0.2}}
\multiput(3, 1)(0.4,0){8}{\line(1,0){0.2}}

\put(6, 3){\oval(1, 5)}

\put(1, 3){\line(1, 1){2}}
\put(1, 3){\line(2, 1){2}}
\put(1, 3){\line(1, 0){2}}
\put(1, 3){\line(1, -1){2}}

\put(3, 5){\line(3, -1){3}}
\put(3, 5){\line(3, -2){3}}
\put(3, 5){\line(3, -4){3}}

\put(3, 4){\line(3, 1){3}}
\put(3, 4){\line(1, -1){3}}
\put(3, 4){\line(3, -1){3}}

\put(3, 3){\line(3, -2){3}}
\put(3, 3){\line(3, 2){3}}
\put(3, 3){\line(3, 1){3}}

\put(3, 1){\line(3, 2){3}}
\put(3, 1){\line(3, 4){3}}
\put(3, 1){\line(1, 1){3}}

\end{picture}
\caption{A $3$-$\gamma_{c}$-critical graph which is non-factor critical}
\label{f:fig001}
\end{center}
\end{figure}
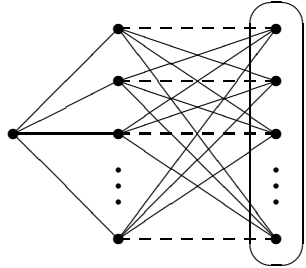

\indent  Chen et al.~\cite{XGC} characterized all $k$-$\gamma_{c}$-critical graphs when $1 \leq k \leq 2$.
\vskip 5 pt

\begin{Theorem}\label{thm chen}
A graph $G$ is $1$-$\gamma_{c}$-critical graph if and only if $G$ is a complete graph. Moreover, a graph $G$ is $2$-$\gamma_{c}$-critical graph if and only if $\overline{G} = \cup^{n}_{i = 1}K_{1, n_{i}}$ where $n \geq 2$ and $n_{i} \geq 1$ for all $1 \leq i \leq n$.
\end{Theorem}
\vskip 5 pt

\noindent We then obtain the following observations as a consequence of theorem \ref{thm chen}.
\vskip 5 pt

\begin{ob}\label{obs 4c}
For $k \in [2]$ every $k$-$\gamma_{c}$-critical graph of odd order with minimum degree $\delta \geq 2$ is factor critical.
\end{ob}

\begin{ob}\label{obs bicri}
For $k \in [2]$ every $k$-$\gamma_{c}$-critical graph of even order with minimum degree $\delta \geq 3$ is bi-critical.
\end{ob}


\vskip 5 pt

\noindent By Theorem \ref{thm chen}, we observe further that a $k$-$\gamma_{c}$-critical graph does not contain a cut vertex when $1 \leq k \leq 2$.
\vskip 5 pt

\begin{ob}\label{obs 0}
Let $G$ be a $k$-$\gamma_{c}$-critical graph with $1 \leq k \leq 2$. Then $G$ has no cut vertex.
\end{ob}

\noindent Further,  Chen et al.\cite{XGC} established fundamental properties of $k$-$\gamma_{c}$-critical graphs for $k \geq 2$.
\vskip 5 pt

\begin{Lemma}\label{lem 1} \cite{XGC}
Let $G$ be a $k$-$\gamma_{c}$-critical graph, $x$ and $y$ a pair of non-adjacent vertices of $G$ and $D_{xy}$ a $\gamma_{c}$-set of $G + xy$. Then
\begin {enumerate}
\item [(1)] $k - 2 \leq |D_{xy}| \leq k - 1$,
\item [(2)] $D_{xy} \cap \{x, y\} \neq \emptyset$,
\item [(3)] if $\{x\} = \{x, y\} \cap D_{xy}$, then $N_{G}(y) \cap D_{xy} = \emptyset$.
\end {enumerate}
\end{Lemma}
\vskip 5 pt

\indent Ananchuen\cite{A} established structures of $k$-$\gamma_{c}$-critical with a cut vertex.
\vskip 5 pt

\begin{Lemma}\label{lem 2}\cite{A}
For $k \geq 3$, let $G$ be a $k$-$\gamma_{c}$-critical graph with a cut vertex $c$ and $D$ a connected dominating set. Then
\begin {enumerate}
\item [(1)] $G - c$ contains exactly two components,
\item [(2)] if $C_{1}$ and $C_{2}$ are the components of $G - c$, then $G[N_{C_{1}}(c)]$ and $G[N_{C_{2}}(c)]$ are complete,
\item [(3)] $c \in D$.
\end {enumerate}
\end{Lemma}
\vskip 5 pt

All the following results of this section were established in \cite{PKNA2}. The first result is the construction of a forbidden subgraph of $k$-$\gamma_{c}$-critical graphs. For a connected graph $G$, let $X, Y, X_{1}$ and $Y_{1}$ be disjoint vertex subsets of $V(G)$. The induced subgraph $G[X \cup X_{1} \cup Y \cup Y_{1}]$ is called a \emph{bad subgraph} if
\begin {enumerate}
\item [(\emph{i})] $x_{1} \succ X \cup X_{1}$ for any vertex $x_{1} \in X_{1}$,
\item [(\emph{ii})] $N[x] \subseteq X \cup X_{1}$ for any vertex $x \in X$,
\item [(\emph{iii})] $y_{1} \succ Y \cup Y_{1}$ for any vertex $y_{1} \in Y_{1}$ and
\item [(\emph{iv})] $N[y] \subseteq Y \cup Y_{1}$ for any vertex $y \in Y$.
\end {enumerate}
\vskip 5 pt

\noindent An example of a bad subgraph is illustrated in Figure \ref{bad}.
\vskip 5 pt

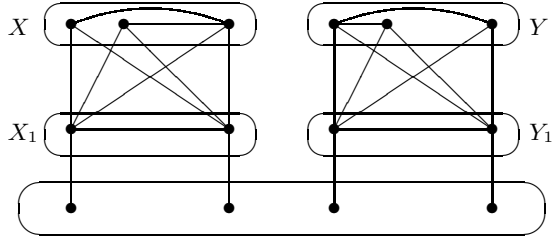
\begin{figure}[htb]
\begin{center}
\setlength{\unitlength}{0.7cm}
\begin{picture}(10, 5)
\put(1, 2){\circle*{0.2}}
\put(1, 4){\circle*{0.2}}
\put(2, 4){\circle*{0.2}}
\put(4, 4){\circle*{0.2}}
\put(4, 2){\circle*{0.2}}
\put(1, 2){\line(0, 1){2}}
\put(1, 2){\line(1, 2){1}}
\put(1, 2){\line(3, 2){3}}
\put(4, 2){\line(0, 1){2}}
\put(4, 2){\line(-1, 1){2}}
\put(4, 2){\line(-3, 2){3}}
\put(6, 2){\circle*{0.2}}
\put(6, 4){\circle*{0.2}}
\put(7, 4){\circle*{0.2}}
\put(9, 4){\circle*{0.2}}
\put(9, 2){\circle*{0.2}}
\put(6, 2){\line(0, 1){2}}
\put(6, 2){\line(1, 2){1}}
\put(6, 2){\line(3, 2){3}}
\put(9, 2){\line(0, 1){2}}
\put(9, 2){\line(-1, 1){2}}
\put(9, 2){\line(-3, 2){3}}

\put(1, 0.5){\circle*{0.2}}
\put(1, 0.5){\line(0, 1){2}}
\put(4, 0.5){\circle*{0.2}}
\put(4, 0.5){\line(0, 1){2}}
\put(6, 0.5){\circle*{0.2}}
\put(6, 0.5){\line(0, 1){2}}
\put(9, 0.5){\circle*{0.2}}
\put(9, 0.5){\line(0, 1){2}}

\qbezier(1, 2)(2.5, 2)(4, 2)
\qbezier(6, 2)(7.5, 2)(9, 2)
\qbezier(1, 4)(2.5, 4.6)(4, 4)
\qbezier(6, 4)(7.5, 4.6)(9, 4)
\put(2, 4){\line(1, 0){2}}
\put(6, 4){\line(1, 0){1}}

\put(5, 0.5){\oval(10, 1)}
\put(2.5, 1.9){\oval(4, 0.8)}
\put(7.5, 1.9){\oval(4, 0.8)}
\put(2.5, 4){\oval(4, 0.8)}
\put(7.5, 4){\oval(4, 0.8)}

\put(-0.2, 1.8){\footnotesize$X_{1}$}
\put(-0.2, 3.8){\footnotesize$X$}
\put(9.7, 1.8){\footnotesize$Y_{1}$}
\put(9.7, 3.8){\footnotesize$Y$}


\end{picture}
\caption{The induced subgraph $G[X \cup X_{1} \cup Y \cup Y_{1}]$}
\label{bad}
\end{center}
\end{figure}
\vskip 20 pt

\noindent The authors showed, in \cite{PKNA2}, that :
\vskip 5 pt

\begin{Lemma}\cite{PKNA2}\label{lem mp0}
For $k \geq 3$, let $G$ be a $k$-$\gamma_{c}$-critical graph. Then $G$ does not contain a bad subgraph.
\end{Lemma}
\vskip 5 pt

\noindent They also provided characterizations of some blocks of $k$-$\gamma_{c}$-critical graphs. Recall that
\vskip 5 pt

\indent $\cC(G)$ is the set of all cut vertices of $G$
\vskip 5 pt

\noindent and, for a block $B$ of $G$,
\vskip 8 pt

\indent $\cC(B) = V(B) \cap \cC(G)$ and $\zeta(G) = |\cC(G)|$.
\vskip 8 pt

\noindent When no ambiguity occur, we write $\cC$ rather than $\cC(G)$. In the same paper, the authors showed further that for a connected graph $G$ and a pair of non-adjacent vertices $x$ and $y$ of $G$, $\cC(G) = \cC(G + xy)$ if $x$ and $y$ are in the same block of $G$.
\vskip 5 pt

\begin{Lemma}\cite{PKNA2}\label{lem Pum1}
For a connected graph $G$, let $B$ be a block of $G$ and $x, y \in V(B)$ such that $xy \notin E(G)$. Then $\cC(G) = \cC(G + xy)$.
\end{Lemma}
\vskip 5 pt

\indent Let $D$ be a $\gamma_{c}$-set of $G$. The followings are the characterization of four classes of end blocks of $k$-$\gamma_{c}$-critical graphs that contains at most $3$ vertices from $D$. For vertices $c, z_{1}$ and $z_{2}$, let
\vskip 5 pt

\indent $\mathcal{B}_{0} = \{c \vee K_{t_{1}} :$ for an integer $t_{1} \geq 1\}$,
\vskip 5 pt

\indent $\mathcal{B}_{1} = \{c \vee K_{t_{2}} \vee z_{1} :$ for an integer $t_{2} \geq 2\}$ and
\vskip 5 pt

\indent $\mathcal{B}_{2, 1} = \{c \vee K_{t_{3}}\vee K_{t_{4}} \vee z_{2} :$ for integers $t_{3}, t_{4} \geq 2\}$.
\vskip 5 pt

\noindent The following is a part of the construction of $\mathcal{B}_{2, 2}$. For integers $l \geq 2, m_{i} \geq 1$ and $r \geq 0$, we let $\mathcal{S} = \cup^{l}_{i = 1}K_{1, m_{i}}$ and $T = \mathcal{S} \cup \overline{K_{r}}$. When $r = 0$, we let $T = \mathcal{S}$. Then, for $1 \leq i \leq l$, let $s^{i}_{0}, s^{i}_{1}, s^{i}_{2}, ..., s^{i}_{m_{i}}$ be the vertices of a star $K_{1, m_{i}}$ which $s^{i}_{0}$ is the center. Further, let $S = \cup^{l}_{i = 1}\{s^{i}_{1}, s^{i}_{2}, ..., s^{i}_{m_{i}}\}$ and $S' = \cup^{l}_{i = 1}\{s^{i}_{0}\}$, moreover, let $S'' = V(\overline{K_{r}})$ if $T = \mathcal{S} \cup \overline{K_{r}}$ and $S'' = \emptyset$ if $T = \mathcal{S}$. Therefore,
\[ \overline{T}      = \left\{\begin{array}{lll}
                     \overline{\mathcal{S}} & or\\
                     \overline{\mathcal{S}} \vee K_{r}.
               \end{array}\right. \]
\noindent That is, $\overline{T}$ can be obtained by removing the edges in the stars of $\mathcal{S}$ from a complete graph on $S \cup S' \cup S''$. Then the blocks in $\mathcal{B}_{2, 2}$ are defined as follows.
\vskip 5 pt

\indent $\mathcal{B}_{2, 2} = \{c \vee \leftidx{_{\overline{T}[S]}}{\overline{T}} :$ for integers $l \geq 2, r \geq 0$ and $m_{i} \geq 1\}$.
\vskip 5 pt

\noindent A graph in this class is illustrated by Figure \ref{blockb}. According to the figure, an \emph{oval} denotes a complete subgraph, \emph{double lines} between subgraphs denote all possible edges between the subgraphs and a \emph{dash line} denotes a removed edge.
\vskip 8 pt

\begin{figure}[htb]
\setlength{\unitlength}{0.8cm}
\begin{center}
\begin{picture}(13, 7)
\put(3, 4){\circle*{0.2}}
\put(7.5, 1.7){\circle*{0.2}}
\put(7.5, 2.7){\circle*{0.2}}
\put(5.5, 4){\oval(1, 4.7)}
\put(7.5, 3.9){\oval(1, 5)}
\put(3.2, 4.1){\line(1, 0){1.6}}
\put(3.2, 3.9){\line(1, 0){1.6}}
\put(7, 3){\line(1, 0){1}}
\put(6.1, 4.9){\line(1, 0){0.8}}
\put(6.1, 4.6){\line(1, 0){0.8}}
\put(6.1, 2.4){\line(1, 0){0.8}}
\put(6.1, 2.1){\line(1, 0){0.8}}
\put(5.5, 6){\circle*{0.2}}
\put(5.5, 5){\circle*{0.2}}
\put(7.5, 6){\circle*{0.2}}
\multiput(7, 6)(-0.7,0){3}{\line(1,0){0.4}}
\multiput(7.5, 6)(-0.7,-0.35){3}{\line(-2,-1){0.5}}
\put(5.45, 5.3){\footnotesize$\vdots$}
\put(7.45, 2){\footnotesize$\vdots$}

\put(5.5, 4){\circle*{0.2}}
\put(5.5, 3){\circle*{0.2}}
\put(7.5, 4){\circle*{0.2}}
\multiput(7, 4)(-0.7,0){3}{\line(1,0){0.4}}
\multiput(7.5, 4)(-0.7,-0.35){3}{\line(-2,-1){0.5}}
\put(5.45, 3.3){\footnotesize$\vdots$}
\put(7.45, 5){\footnotesize$\vdots$}

\put(2.8, 3.5){\footnotesize$c$}
\put(5.45, 1.2){\footnotesize$S$}
\put(7.45, 0.9){\footnotesize$S''$}
\put(8.3, 5){\footnotesize$S'$}
\end{picture}
\caption{A block $B$ in the class $\mathcal{B}_{2, 2}$}
\label{blockb}
\end{center}
\end{figure}
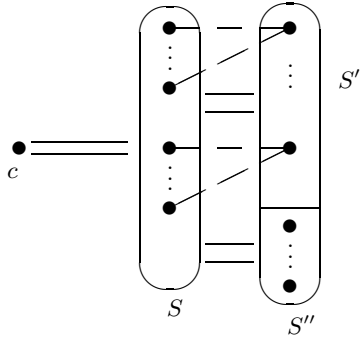

\vskip 20 pt

\indent For a block $B \in \mathcal{B}_{0} \cup \mathcal{B}_{1} \cup \mathcal{B}_{2, 1} \cup \mathcal{B}_{2, 2}$, the vertex $c$ is called the \emph{head} of $B$. The following are the characterizations of an end block $B$ such that $|D \cap V(B)| \leq 3$.

\begin{Lemma}\cite{PKNA2}\label{lem pmp0}
Let $G$ be a $k$-$\gamma_{c}$-critical graph with a $\gamma_{c}$-set $D$ and let $B$ be an end block of $G$. If $|(D \cap V(B)) - \cC| = 0$, then $B \in \mathcal{B}_{0}$.
\end{Lemma}
\vskip 5 pt

\begin{Lemma}\cite{PKNA2}\label{lem pmp1}
Let $G$ be a $k$-$\gamma_{c}$-critical graph with a $\gamma_{c}$-set $D$ and let $B$ be an end block of $G$. If $|(D \cap V(B)) - \cC| = 1$, then $B \in \mathcal{B}_{1}$.
\end{Lemma}
\vskip 5 pt

\begin{Lemma}\cite{PKNA2}\label{lem pmp2}
Let $G$ be a $k$-$\gamma_{c}$-critical graph with a $\gamma_{c}$-set $D$ and let $B$ be an end block of $G$. Suppose that $|(D \cap V(B)) - \cC| = 2$. Then $B \in \mathcal{B}_{2, 1} \cup \mathcal{B}_{2, 2}$.
\end{Lemma}
\vskip 5 pt

\begin{Lemma}\cite{PKNA2}\label{lem pum}
For $k \geq 3$, let $G$ be a $k$-$\gamma_{c}$-critical graph. Then $G$ contains at most one end block $B$ such that $B \in \mathcal{B}_{0} \cup \mathcal{B}_{1}$.
\end{Lemma}
\vskip 5 pt

\noindent Finally, we conclude this section by the following two lemmas which are structures of blocks of $k$-$\gamma_{c}$-critical graphs. Recall from Section \ref{main} that
\vskip 8 pt

\indent $\mathcal{Z}(k, \zeta) :$ the class of $k$-$\gamma_{c}$-critical graphs containing $\zeta$ cut vertices.
\vskip 8 pt

\begin{Lemma}\cite{PKNA2}\label{lem b}
Let $G \in \mathcal{Z}(k, \zeta)$ where $\zeta \in \{k - 3, k - 2\}$. Then $G$ has only two end blocks and another blocks contain exactly two cut vertices.
\end{Lemma}
\vskip 5 pt

\begin{Lemma}\cite{PKNA2}\label{lem c}
Let $G \in \mathcal{Z}(k, \zeta)$ with a $\gamma_{c}$-set $D$ where $\zeta \in \{k - 3, k - 2\}$ and $B$ be a block of $G$ containing two cut vertices $c$ and $c'$. If $(D \cap V(B)) - \cC = \emptyset$. Then $B = cc'$.
\end{Lemma}
\vskip 8 pt

\indent In \cite{HAK}, the authors established a construction of the class $\mathcal{P}(k)$ which, for any graph $G \in \mathcal{P}(k)$ and integer $l \geq 1$, there exists a $(k + l)$-$\gamma_{c}$-critical graph contains $G$ as an induced subgraph. Recall that, for a subgraph $H$ of $G$, $H$ is maximal complete subgraph of $G$ if for any complete subgraph $H'$ of $G$ such that $V(H) \subseteq V(H')$, we have $V(H) = V(H')$.
\vskip 5 pt

\noindent \textbf{The class $\mathcal{P}(k)$}\\
\indent A $k$-$\gamma_{c}$-critical graph $G$ is in this class if there exists a maximal complete subgraph $H$ of order at least two of $G$ satisfies the following properties :
\begin{enumerate}
\item [(\emph{i})] for any vertex $x$ of $G$, there exists a $\gamma_{c}$-set $D$ of $G$ such that $x \in D$ and $D \cap V(H) \neq \emptyset$ and
\item [(\emph{ii})] for any non-adjacent vertices $x$ and $y$ of $G$, $D'_{xy} \cap V(H) \neq \emptyset$ where $D'_{xy}$ is a connected dominating set of $G + xy$ such that $|D'_{xy}| < k$ ($D'_{xy}$ need not be a $\gamma_{c}$-set of $G + xy$).
\end{enumerate}
\vskip 5 pt

\noindent We next give a construction of a $(k + l)$-$\gamma_{c}$-critical graph containing $G$ in the class $\mathcal{P}(k)$ as an induced subgraph. Let $H$ be a maximal complete subgraph of $G$ having properties $(i)$ and $(ii)$. The graph $G(n_{1}, n_{2}, ..., n_{l})$ can be constructed from a vertex $x_{0}$, $l$ copies of completes graph $K_{n_{1}}, K_{n_{2}}, ..., K_{n_{l- 1}}$ which $n_{i} \geq 1$ for $1 \leq i \leq l$ and the graph $G$ by adding edges according to the join operations :
\vskip 5 pt

\indent $x_{0} \vee K_{n_{1}} \vee K_{n_{2}} \vee ... \vee K_{n_{l}} \vee \leftidx{_H}{G}$.
\vskip 5 pt

\noindent The graph is illustrated by Figure \ref{join}.

\vskip 20 pt

\begin{figure}[htb]
\begin{center}
\setlength{\unitlength}{0.7cm}
\begin{picture}(20, 5)
\put(3.8, 3){\circle*{0.2}}
\put(3.8, 3){\line(1, 1){0.7}}
\put(3.8, 3){\line(1, -1){0.7}}
\put(9, 3){\oval(1, 2)}
\put(9.6, 2.5){\line(1, 0){0.7}}
\put(9.6, 3.5){\line(1, 0){0.7}}
\put(9.6, 3){\line(1, 0){0.7}}
\put(5, 3){\oval(1, 2)}
\put(5.6, 2.5){\line(1, 0){0.7}}
\put(5.6, 3.5){\line(1, 0){0.7}}
\put(5.6, 3){\line(1, 0){0.7}}
\put(7, 3){\oval(1, 2)}
\put(11, 3){\oval(1, 2)}
\put(12, 3){\oval(3, 3)}
\put(3.7, 2.2){\footnotesize$x_{0}$}
\put(5, 1.5){\footnotesize$K_{n_{l}}$}
\put(6.8, 1.5){\footnotesize$K_{n_{2}}$}
\put(8.8, 1.5){\footnotesize$K_{n_{l}}$}
\put(7.8, 2.8){\footnotesize$...$}
\put(11.3, 1){\footnotesize$G \in \mathcal{P}(k)$}
\put(10.8, 2.5){\footnotesize$H$}
\end{picture}
\caption{The graph $G(n_{1}, n_{2}, ..., n_{l})$}
\label{join}
\end{center}
\end{figure}
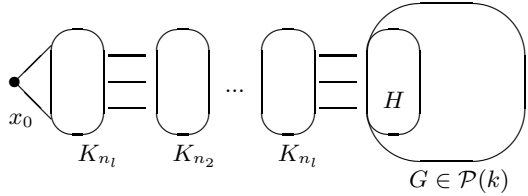
\vskip 10 pt

\noindent Thus, they proved that

\begin{Theorem}\label{thm gl}\cite{HAK}
For an integer $k \geq 3$, let $G \in \mathcal{P}(k)$. Then $G(n_{1}, n_{2}, ..., n_{l})$ is a $(k + l)$-$\gamma_{c}$-critical graph for all $l \geq 1$.
\end{Theorem}

\section{Proofs}\label{proof}

\subsection{Connected Dominating Set of Blocks}\label{block}
\noindent Let
\vskip 5 pt

\indent $\mathfrak{B}(G)$ be the family of all blocks of $G$.
\vskip 8 pt

\noindent When no ambiguity can occur, we use $\mathfrak{B}$ to denote $\mathfrak{B}(G)$. For a $k$-$\gamma_{c}$-critical graph $G$ with a cut vertex, let $B$ be an end block of $G$ containing non-adjacent vertices $x$ and $y$. Clearly, $V(B + xy) = V(B)$. Let $D$ be a $\gamma_{c}$-set of $G$.
\vskip 5 pt

\begin{Lemma}\label{cor Pum1}
Let $B$ be a block of $G$ and $x, y \in V(B)$ such that $xy \notin E(G)$. Then $D \cap \cC = D_{xy} \cap \cC$, in particular, $D \cap \cC(B') = D_{xy} \cap \cC(B')$ for all $B' \in \mathfrak{B}(G + xy)$.
\end{Lemma}

\proof
We first show that $D \cap \cC \subseteq D_{xy} \cap \cC$. Let $c \in D \cap \cC$. By Lemma \ref{lem Pum1}, $c \in \cC(G + xy)$. By the connectedness of $(G + xy)[D_{xy}]$, $c \in D_{xy}$. Thus $D \cap \cC \subseteq D_{xy} \cap \cC$. We now show that $D_{xy} \cap \cC \subseteq D \cap \cC$. Let $c \in D_{xy} \cap \cC$. That is $c \in \cC$. Lemma \ref{lem 2}(3) yields that $c \in D$. So $c \in D \cap \cC$ and thus, $D_{xy} \cap \cC \subseteq D \cap \cC$, as required.
\vskip 5 pt

\indent In view of Lemma \ref{lem Pum1}, $V(B') \cap \cC(G + xy) = V(B') \cap \cC$ for all $B' \in \mathfrak{B}(G + xy)$. Because $D \cap \cC = D_{xy} \cap \cC$, it follows that
\begin{center}
$D \cap \cC(B') = D \cap \cC \cap V(B') = D_{xy} \cap \cC \cap V(B') = D_{xy} \cap \cC(B')$.
\end{center}
This completes the proof.
\qed

\vskip 5 pt

It is worth noting that, in \cite{PKNA2}, the similar result as Lemma \ref{cor Pum1} was proved but focused only end blocks. So, our result in Lemma \ref{cor Pum1} is more general. For non-adjacent vertices $x$ and $y$ of the block $B$, the following lemma gives the number of vertices of a $\gamma_{c}$-set of $G + xy$ in $B$.
\vskip 5 pt

\begin{Lemma}\label{lem Pum2}
For all $x, y \in V(B)$ such that $xy \notin E(G)$, $|D_{xy} \cap V(B)| < |D \cap V(B)|$.
\end{Lemma}

\proof
We first establish the following claim.
\vskip 5 pt

\noindent \textbf{Claim :} For all block $B'$ which is not $B$, $|D \cap V(B')| \leq |D_{xy} \cap V(B')|$.\\
\indent Suppose to the contrary that $|D \cap V(B')| > |D_{xy} \cap V(B')|$. Lemma \ref{cor Pum1} gives that $D \cap \cC(B') = D_{xy} \cap \cC(B')$. Because $x, y \notin V(B')$, $G[(D - V(B')) \cup (D_{xy} \cap V(B'))]$ is connected. Moreover, $(D - V(B')) \cup (D_{xy} \cap V(B')) \succ_{c} G$. This implies that
\begin{align}
k = |D| = |(D - V(B')) \cup (D \cap V(B'))| & = |D - V(B')| + |D \cap V(B')|\notag\\
                                             & > |D - V(B')| + |D_{xy} \cap V(B')|\notag\\
                                             & = |(D - V(B')) \cup (D_{xy} \cap V(B'))|,\notag
\end{align}
\noindent contradicting the minimality of $D$. Thus establishing the claim.
\vskip 5 pt

\indent We now prove Lemma \ref{lem Pum2}. Suppose to the contrary that $|D_{xy} \cap V(B)| \geq |D \cap V(B)|$. Lemma \ref{cor Pum1} yields that $D \cap \cC = D_{xy} \cap \cC$. Clearly $D = \cup_{\tilde{B} \in \mathfrak{B}}(D \cap V(\tilde{B}))$ and $D_{xy} = \cup_{\tilde{B} \in \mathfrak{B}(G + xy)}(D_{xy} \cap V(\tilde{B}))$. Lemma \ref{lem 2}(1) yields, further, that each cut vertex is contained in exactly two blocks. Thus each cut vertex is counted twice in $\sum_{\tilde{B} \in \mathfrak{B}}|(D \cap V(\tilde{B}))|$ and $\sum_{\tilde{B} \in \mathfrak{B}(G + xy)}|(D_{xy} \cap V(\tilde{B}))|$. Therefore, $|D| = \Sigma_{\tilde{B} \in \mathfrak{B}}|D \cap V(\tilde{B})| - |\cC|$ and $\Sigma_{\tilde{B} \in \mathfrak{B}(G + xy)}|D_{xy} \cap V(\tilde{B})| - |\cC(G + xy)| = |D_{xy}|$. We note by Lemma \ref{lem Pum1} that $|\cC| = |\cC(G + xy)|$. By the claim and the assumption that $|D_{xy} \cap V(B)| \geq |D \cap V(B)|$, we have
\begin{align}
k = |D| & =     \Sigma_{\tilde{B} \in \mathfrak{B}}|D \cap V(\tilde{B})| - |\cC|\notag\\
        & =     |D \cap V(B)| + \Sigma_{\tilde{B} \in \mathfrak{B} - \{B\}}|D \cap V(\tilde{B})| - |\cC|\notag\\
        & \leq  |D_{xy} \cap V(B)| + \Sigma_{\tilde{B} \in \mathfrak{B} - \{B\}}|D \cap V(\tilde{B})| - |\cC|~~(by~the~assumption)\notag\\
        & \leq  |D_{xy} \cap V(B)| + \Sigma_{\tilde{B} \in \mathfrak{B}(G + xy) - \{B\}}|D_{xy} \cap V(\tilde{B})| - |\cC|~~(by~the~claim)\notag\\
        & =     \Sigma_{\tilde{B} \in \mathfrak{B}(G + xy)}|D_{xy} \cap V(\tilde{B})| - |\cC(G + xy)| =     |D_{xy}|.\notag
\end{align}
\noindent This contradicts Lemma \ref{lem 1}(1). Thus $|D_{xy} \cap V(B)| < |D \cap V(B)|$ and this completes the proof.
\qed

\vskip 5 pt

\begin{Corollary}\label{cor Pum3}
For all block $B$ of $G$ and $x, y \in V(B)$ such that $xy \notin E(G)$, $|(D_{xy} \cap V(B)) - \cC| < |(D \cap V(B)) - \cC|$.
\end{Corollary}
\proof
In view of Lemma \ref{cor Pum1}, $D \cap V(B) \cap \cC = D_{xy} \cap V(B) \cap \cC$. Lemma \ref{lem Pum2} then implies that
\begin{align}
|(D_{xy} \cap V(B)) - \cC| & = |D_{xy} \cap V(B)| - |D_{xy} \cap V(B) \cap \cC|\notag\\
                                    & < |D \cap V(B)| - |D \cap V(B) \cap \cC|\notag\\
                                    & = |(D \cap V(B)) - \cC|\notag
\end{align}

\noindent and this completes the proof.
\qed

\vskip 5 pt

\subsection{The $k$-$\gamma_{c}$-Critical Graphs with $\zeta(G) = k - 3$}\label{k-3}
\noindent In this subsection, we characterize $k$-$\gamma_{c}$-critical graphs with $k - 3$ cut vertices. We recall the classes $\mathcal{G}_{1}(i_{1}, i_{2}, ..., i_{k - 3})$ and $\mathcal{G}_{2}(k)$ from Section \ref{main}. First, we will prove that all graphs in these two classes are $k$-$\gamma_{c}$-critical with $k - 3$ cut vertices.
\vskip 5 pt

\begin{Lemma}\label{lem mp5}
Let $G$ be a graph in the class $\mathcal{G}_{1}(i_{1}, i_{2}, ..., i_{k - 3})$, then $G$ is a $k$-$\gamma_{c}$-critical graph with $k - 3$ cut vertices.
\end{Lemma}

\proof
Clearly $G$ has $c_{1}, c_{2}, ..., c_{k - 4}$ and $c$ as the $k - 3$ cut vertices. Observe that, for any $\textbf{\emph{i}} = (i_{1}, i_{2}, ..., ,i_{k - 3})$, a graph $G \in \mathcal{G}_{1}\textbf{\emph{i}}$ has the path $P = c_{0}, c_{1}, ..., c_{l - 1}, a, c_{l}, ..., c_{k - 4}, c$ from $c_{0}$ to $c$ where $a \in V(K_{n_{l}})$. To prove all cases of $\textbf{\emph{i}}$, we may relabel the path $P$ to be $x_{1}, ..., x_{k - 1}$. Hence, $c_{0} = x_{1}, c_{1} = x_{2}, ...c_{k - 4} = x_{k - 2}$ and $c = x_{k - 1}$. We see that $\{x_{2}, x_{3}, ..., x_{k - 2}, x_{k - 1}, s^{1}_{1}, s^{2}_{0}\} \succ_{c} G$ where $s^{1}_{1}, s^{2}_{0}$ are vertices in $B \in \mathcal{B}_{2, 2}$ defined in Section \ref{main}. Therefore $\gamma_{c}(G) \leq k$.
\vskip 5 pt

\indent Let $D$ be a $\gamma_{c}$-set of $G$. If $x_{1} \notin D$, then, to dominate $x_{1}$, $x_{2} \in D$ . If $x_{1} \in D$, then $x_{2} \in D$ since $G[D]$ is connected. In both cases, $x_{2} \in D$.
\vskip 5 pt

\indent Suppose that $D \cap S'' = \emptyset$. Because $B \in \mathcal{B}_{2, 2}$, to dominate $B$, $|D \cap (S \cup S')| \geq 2$. By the connectedness of $G[D]$, we have $\{x_{3}, x_{4}, ..., x_{k - 2}, x_{k - 1}\} \subseteq D$. Thus $\gamma_{c}(G) \geq k$ implying that $\gamma_{c}(G) = k$. Hence, suppose that $D \cap S'' \neq \emptyset$. Since $x_{2} \in D$ and $G[D]$ is connected, it follows that $\{x_{3}, x_{4}, ..., x_{k - 2},$ $x_{k - 1}, y\} \subseteq D$ where $y \in D \cap S$. Thus $\gamma_{c}(G) = |D| \geq k$ implying that $\gamma_{c}(G) = k$.
\vskip 5 pt

\indent Now, we will establish the criticality. Let $u$ and $v$ be a pair of non-adjacent vertices of $G$ and let $S_{1} = S \cup S' \cup S''$. We first assume that $|\{u, v\} \cap S_{1}| = 0$. Therefore $\{u, v\} \subseteq \{x_{1}, x_{2}, ..., x_{k - 2},$ $x_{k - 1}\}$. Thus $u = x_{i}$ and $v = x_{j}$ for some $i, j \in \{1, 2, ..., k - 1\}$. Without loss of generality let $i < j$. Clearly $i + 2 \leq j$. We see that
\begin{center}
$\{x_{2}, x_{3}, ..., x_{i}, x_{i + 2}, x_{i + 3}, ..., x_{j}, ..., x_{k - 1}, s^{1}_{1}, s^{2}_{0}\} \succ_{c} G + uv$.
\end{center}
\noindent So $\gamma_{c}(G + uv) \leq k - 1$.
\vskip 5 pt

\indent Hence, we assume that $|\{u, v\} \cap S_{1}| = 1$. If $\{u, v\} = \{x_{k - 1}, s\}$ for some $s \in S_{1}$, then $s \notin S$ and, clearly, $\{x_{k - 1}, s\} \succ S_{1}$. Thus $\{x_{2}, x_{3}, ..., x_{k - 1}, s_{1}\} \succ_{c} G + uv$. Therefore $\gamma_{c}(G + uv) \leq k - 1$. Let $v \in S_{1}$. Since $|\mathcal{S}| \geq 2$, there exists $v' \in S - \{v\}$ such that $\{v, v'\} \succ_{c} S_{1}$. Suppose that $u \in \{x_{2}, x_{3}, ..., x_{k - 2}\}$. Thus $\{x_{2}, x_{3}, ..., u, ..., x_{k - 2}, v, v'\} \succ_{c} G + uv$. Hence $\gamma_{c}(G + uv) \leq k - 1$. If $u = x_{1}$, then $\{x_{3}, x_{4}, ..., x_{k - 2}, v, v'\} \succ_{c} G + uv$ implying that $\gamma_{c}(G + uv) \leq k - 1$.
\vskip 5 pt

\indent Finally, we assume that $\{u, v\} \subseteq S_{1}$. Thus $\{u, v\} = \{s^{i}_{j}, s^{i}_{0}\}$ for some $i \in \{1, 2, ..., |\mathcal{S}|\}$ and $j \in \{1, 2, ..., m_{i}\}$. Clearly $\{x_{2}, x_{3}, ..., x_{k - 1}, s^{i}_{j}\} \succ_{c} G + uv$ and $\gamma_{c}(G + uv) \leq k - 1$. Thus $G$ is a $k$-$\gamma_{c}$-critical graph and this completes the proof.
\qed

\vskip 8 pt

\begin{Lemma}\label{lem d}
Let $G$ be a graph in the class $\mathcal{G}_{2}(k)$. Then $G$ is a $k$-$\gamma_{c}$-critical graph with $k - 3$ cut vertices.
\end{Lemma}

\proof
Choose $v \in V(\check{B})$. By $(1)$ of $\mathcal{B}_{3}$, there exists a $\gamma_{c}$-set $D_{v}$ such that $|D_{v}| = 3$ and $D_{v} \cap A \neq \emptyset$. Thus $\{c_{1}, c_{2}, ..., c_{k - 4}, c_{k - 3}\} \cup D_{v} \succ_{c} G$. Therefore $\gamma_{c}(G) \leq k$.
\vskip 5 pt

\indent Let $D$ be a $\gamma_{c}$-set of $G$. As $c_{1}, c_{2}, ..., c_{k - 3}$ are cut vertices, by Lemma \ref{lem 2}(3), $c_{1}, c_{2}, ..., c_{k - 3} \in D$. Let $v \in V(\check{B}) \cap D$. Observe that $V(\check{B}) \cap D$ is a connected dominating set of $B$ containing $v$. By $(1)$ and the minimality of $D_{v}$, $|V(\check{B}) \cap D| \geq |D_{v}| = 3$. So $\gamma_{c}(G) \geq k$ and this implies that $\gamma_{c}(G) = k$.
\vskip 5 pt

\indent We will prove the criticality. Let $u$ and $v$ be non-adjacent vertices of $G$. Suppose first that $c_{0} \in \{u, v\}$, $c_{0} = u$ say. If $v \in \{c_{2}, c_{3}, ..., c_{k - 3}\}$, then $\{c_{2}, c_{3}, ..., c_{k - 3}\} \cup D_{v} \succ_{c} G + uv$. If $v \in V(\check{B})$, then, by $(1)$, there exists a $\gamma_{c}$-set $D_{v}$ of size $3$ of $B$ such that $v \in D_{v}$ and $A \cap D_{v} \neq \emptyset$. So $\{c_{2}, c_{3}, ..., c_{k - 3}\} \cup D_{v} \succ_{c} G + uv$. These imply that $\gamma_{c}(G + uv) < \gamma_{c}(G)$.
\vskip 5 pt

\indent We then suppose that $c_{0} \notin \{u, v\}$. If $\{u, v\} \subseteq \{c_{1}, c_{2}, ..., c_{k - 3}\}$, then there exists $i$ and $j$ such that $c_{i} = u$ and $c_{j} = v$. Without loss of generality let $i < j$. Clearly, $i + 2 \leq j$. So $\{c_{1}, c_{2}, ..., c_{i}, c_{i + 2}, c_{i + 3},$ $..., c_{k - 3}\} \cup D_{v} \succ_{c} G + uv$. For the case when $|\{u, v\} \cap \{c_{1}, c_{2}, ..., c_{k - 3}\}| = 1$, we have $\{c_{1}, ..., c_{i}, ..., c_{k - 4}\} \cup D_{v} \succ_{c} G + uv$. Finally, if $\{u, v\} \subseteq V(B)$, then, by $(2)$, there exists a $\gamma_{c}$-set $D^{B}_{uv}$ such that $D^{B}_{uv} \cap \{u, v\} \neq \emptyset$, $|D^{B}_{uv}| = 2$ and $D^{B}_{uv} \cap A \neq \emptyset$. Thus $\{c_{1}, c_{2}, ..., c_{k - 3}\} \cup D^{B}_{uv} \succ_{c} G + uv$. This implies that $\gamma_{c}(G + uv) < \gamma_{c}(G)$. Clearly, $c_{1}, c_{2}, .., c_{k - 3}$ are the $k - 3$ cut vertices of $G$. This completes the proof.
\qed

\vskip 5 pt

\indent In the following, we let $G \in \mathcal{Z}(k, k - 3)$ having a $\gamma_{c}$-set $D$. In view of Lemma \ref{lem b}, $G$ has only two end blocks and another blocks contain two cut vertices. Thus, we let $B_{1}$ and $B_{k - 2}$ be the two end blocks and another blocks $B_{2}, B_{3}, ..., B_{k - 3}$ contain two cut vertices. Without loss of generality let $c_{1} \in V(B_{1}), c_{k - 3} \in V(B_{k - 2})$ and $c_{i - 1}, c_{i} \in V(B_{i})$ for $2 \leq i \leq k - 3$. Moreover, let $C_{i} = V(B_{i}) - \cC$ for all $1 \leq i \leq k - 2$. Let $D'$ be a $\gamma_{c}$-set of $G$ such that $D' \neq D$, by the minimality of $k$, we have $|V(B_{i}) \cap D| = |V(B_{i}) \cap D'|$ for all $i$. Thus, we can let
\vskip 10 pt

\indent $\mathcal{H}(b_{1}, b_{2}, b_{3}, ..., b_{k - 2}) :$ the class of a graph $G \in \mathcal{Z}(k, k - 3)$ such that\\
\indent $\quad \quad \quad \quad \quad \quad \quad \quad \quad \quad$ $|V(C_{i}) \cap D| = b_{i}$ for $1 \leq i \leq k - 2$.
\vskip 10 pt

\begin{Lemma}\label{lem e}
For a $\gamma_{c}$-set of $G$, either $|V(C_{1}) \cap D| \geq 2$ or $|V(C_{k - 2}) \cap D| \geq 2$.
\end{Lemma}

\proof
Suppose to the contrary that $|V(C_{1}) \cap D| \leq 1$ and $|V(C_{k - 2}) \cap D| \leq 1$. Lemmas \ref{lem pmp0} and \ref{lem pmp1} imply that $B_{1}, B_{k - 2} \in \mathcal{B}_{0} \cup \mathcal{B}_{1}$. This contradicts Lemma \ref{lem pum}. Thus either $|V(C_{1}) \cap D| \geq 2$ or $|V(C_{k - 2}) \cap D| \geq 2$ and this completes the proof.
\qed

\vskip 5 pt

By Lemma \ref{lem e}, we may suppose without loss of generality that $|V(C_{k - 2}) \cap D| \geq |V(C_{1}) \cap D|$.

\begin{Lemma}\label{lem e2}
$\mathcal{Z}(k, k - 3) = \mathcal{H}(0, 0, 0, ..., 3) \cup \mathcal{H}(b_{1}, b_{2}, ..., b_{k - 3}, 2)$ where $b_{i} = 1$ for some $1 \leq i \leq k - 3$ and $b_{j} = 0$ for all $1 \leq j \neq i \leq k - 3$.
\end{Lemma}

\proof
By the definition, $\mathcal{H}(0, 0, 0, ..., 3) \cup \mathcal{H}(b_{1}, b_{2}, ..., b_{k - 3}, 2) \subseteq \mathcal{Z}(k, k - 3)$.
\vskip 5 pt

\indent Conversely, let $G \in \mathcal{Z}(k, k - 3)$. Thus, by Lemma \ref{lem b}, $G$ has only two end blocks $B_{1}$ and $B_{k - 2}$ and another blocks $B_{2}, B_{3}, ..., B_{k - 2}$ contain two cut vertices. Moreover, $c_{1} \in V(B_{1}), c_{k - 3} \in V(B_{k - 1})$ and $c_{i - 1}, c_{i} \in V(B_{i})$ for $2 \leq i \leq k - 3$. In view of Lemma \ref{lem 2}(3), $c_{1}, c_{2}, ..., c_{k - 3}$ $\in D$. So $|D - \cC| = 3$. Thus, $|V(C_{i}) \cap D| \leq 3$ for all $i \in \{1, k - 2\}$. Recall that $|V(C_{k - 2}) \cap D| \geq |V(C_{1}) \cap D|$. Lemma \ref{lem e} implies that either $|V(C_{k - 2}) \cap D| = 3$ or $|V(C_{k - 2}) \cap D| = 2$. That is $G \in \mathcal{H}(0, 0, 0, ..., 3) \cup \mathcal{H}(b_{1}, b_{2}, ..., b_{k - 3}, 2)$. This completes the proof.
\qed

\indent By Lemma \ref{lem e2}, to characterize a graph $G$ in the class $\mathcal{Z}(k, k - 3)$, it suffices to consider when $G$ is either in $\mathcal{H}(0, 0, 0, ..., 3)$ or $\mathcal{H}(b_{1}, b_{2}, ..., b_{k - 3}, 2)$. We first consider the case when $G \in \mathcal{H}(b_{1}, b_{2}, ...,$ $b_{k - 3}, 2)$. Let $c_{i}$ and $c_{i + 1}$ be vertices and $K_{n_{i}}$ a copy of a complete graph.
\vskip 5 pt

\begin{Lemma}\label{lem f}
Let $G \in \mathcal{H}(b_{1}, b_{2}, ..., b_{k - 3}, 2)$ with a block $B_{i}$ containing two cut vertices $c_{i - 1}$ and $c_{i}$ and $b_{i} = 1$. Then $B_{i} = c_{i - 1} \vee K_{n_{i}} \vee c_{i}$ where $n_{i} \geq 2$.
\end{Lemma}

\proof
As $G \in \mathcal{H}(b_{1}, b_{2}, ..., b_{k - 3}, 2)$ and $b_{i} = 1$ for some $1 \leq i \leq k - 3$, we must have $b_{j} = 0$ for all $1 \leq j \neq i \leq k - 3$. Because $B_{i}$ contains two cut vertices, $i > 1$. Therefore, $b_{1} = 0$. Lemma \ref{lem pmp0} then implies that
\vskip 5 pt

\indent $B_{1} = K_{n_{1}} \vee c_{1}$.
\vskip 5 pt

\noindent Let $B' = B_{i} - c_{i - 1} - c_{i}$. We first show that $B'$ is complete. Let $x$ and $y$ be non-adjacent vertices of $B'$. Consider $G + xy$. Lemma \ref{lem 1}(2) implies that $|D_{xy} \cap \{x, y\}| \geq 1$. As $x, y \in V(B')$, we must have $|D_{xy} \cap (V(B_{i}) - \cC)| \geq 1 = b_{i} = |D \cap (V(B_{i}) - \cC)|$ contradicting Corollary \ref{cor Pum3}. So $B'$ is complete.
\vskip 5 pt

\indent We will show that $c_{i - 1}c_{i} \notin E(G)$. Hence, we may assume to the contrary that $c_{i - 1}c_{i} \in E(G)$.  We let
\vskip 5 pt

\indent $X_{1} = N_{B_{i}}(\{c_{i - 1}, c_{i}\})$ and
\vskip 5 pt

\indent $X = V(B') - X_{1}$.
\vskip 5 pt

\noindent Since $|D \cap (V(B_{i}) - \{c_{i - 1}, c_{i}\})| = 1$, it follows that $X \neq \emptyset$. Because $B'$ is complete, $G[X_{1} \cup X]$ is complete. In fact, $X_{1}$ and $X$ satisfy $(i)$ and $(ii)$ of bad subgraphs. We then let
\vskip 5 pt

\indent $Y_{1} = \{c_{1}\}$ and
\vskip 5 pt

\indent $Y = V(K_{n_{1}})$.
\vskip 5 pt

\noindent Thus $G$ has $X, X, Y$ and $Y_{1}$ as a bad subgraph. This contradicts Lemma \ref{lem mp0}. Hence, $c_{i - 1}c_{i} \notin E(G)$.
\vskip 5 pt

\indent We finally show that $N_{B_{i}}(c_{i}) = N_{B_{i}}(c_{i - 1}) = V(B')$. We may assume to the contrary that there exists a vertex $u$ of $B'$ which is not adjacent to $c_{i - 1}$. Consider $G + uc_{i - 1}$. Corollary \ref{cor Pum3} gives that $|D_{uc_{i - 1}} \cap V(B')| = 0$. Lemma \ref{lem Pum1} gives further that $\{c_{i - 1}, c_{i}\} \subseteq D_{uc_{i - 1}}$. Since $c_{i - 1}c_{i} \notin E(G)$, it follows that $(G + uc_{i - 1})[D_{uc_{i - 1}}]$ is not connected, a contradiction. Hence, $N_{B_{i}}(c_{i - 1}) = V(B')$ and, similarly, $N_{B_{i}}(c_{i}) = V(B')$. Since $B_{i}$ is a block, $n_{i} \geq 2$ and this completes the proof.
\qed
\vskip 5 pt

\indent We will prove the following two theorems, both of which give main contribution to the proof of our first main theorem, Theorem \ref{thm k-3}.

\begin{Theorem}\label{thm giii}
For $k \geq 4$, $\mathcal{H}(b_{1}, b_{2}, ..., b_{k - 3}, 2) = \mathcal{G}_{1}(i_{1}, i_{2}, ..., i_{k - 3})$ where $b_{j} = i_{j}$ for all $1 \leq j \leq k - 3$.
\end{Theorem}

\proof
Let $b_{j} = i_{j}$ for all $1 \leq j \leq k - 3$. In views of Lemma \ref{lem mp5}, we have that $\mathcal{G}_{1}(i_{1}, i_{2}, ..., i_{k - 3}) \subseteq \mathcal{H}(b_{1}, b_{2}, ..., b_{k - 3}, 2)$. Thus, it suffices to show that $\mathcal{H}(b_{1}, b_{2}, ...,$ $b_{k - 3}, 2) \subseteq \mathcal{G}_{1}(i_{1}, i_{2}, ..., i_{k - 3})$.
\vskip 5 pt

\indent We will show that $B_{k - 2} \in \mathcal{B}_{2, 2}$. Clearly, $b_{1}$ is either $0$ or $1$. If $b_{1} = 0$, then Lemma \ref{lem pmp0} implies that $B_{1} = K_{n_{1}} \vee c_{1}$. But if $b_{1} = 1$, then Lemma \ref{lem pmp1} implies that $B_{1} = c_{0} \vee K_{n_{1}} \vee c_{1}$. Thus, we let
\[ X_{1}  = \left\{\begin{array}{ll}
                     \{c_{1}\}       & if~~b_{1} = 0~~and\\
                     V(K_{n_{1}})    & if~~b_{1} = 1.
               \end{array}\right. \]
\noindent and
\[ X  = \left\{\begin{array}{ll}
                     V(K_{n_{1}}) & if~~b_{1} = 0~~and\\
                     \{c_{0}\}        & if~~b_{1} = 1.
               \end{array}\right. \]
\noindent Since $b_{k - 2} = 2$, by Lemma \ref{lem pmp2}, $B_{k - 2} \in \mathcal{B}_{2, 1} \cup \mathcal{B}_{2, 2}$. If $B_{k - 2} \in \mathcal{B}_{2, 1}$, then $B_{k - 2} = c_{k - 3} \vee K_{n_{1}} \vee K_{n_{2}} \vee z_{2}$ where $z_{2}$ is given at the definition of $\mathcal{B}_{2, 1}$. We then let
\vskip 5 pt

\indent $Y_{1} = V(K_{n_{2}})$ and $Y = \{z_{2}\}$
\vskip 5 pt

\noindent Clearly, $G$ has a bad subgraph, contradicting Lemma \ref{lem mp0}. Thus $B_{k - 2} \in \mathcal{B}_{2, 2}$.
\vskip 5 pt

\indent We now consider the case when $b_{1} = 1$. Thus $b_{2} = b_{3} = ... = b_{k - 3} = 0$. By Lemma \ref{lem pmp1}, $B_{1} \in \mathcal{B}_{1}$ implying that $B_{1} = c_{0} \vee K_{n_{1}} \vee c_{1}$. Further, Lemma \ref{lem c} implies also that $B_{i} = c_{i - 1}c_{i}$ for $2 \leq i \leq k - 3$. Thus $G \in \mathcal{G}_{1}(1, 0, ..., 0)$.
\vskip 5 pt

\indent We finally consider the case when $b_{1} = 0$. Thus $b_{j} = 1$ for some $2 \leq j \leq k - 3$ and $b_{j} = 0$ for $2 \leq j' \neq j \leq k - 3$. Similarly, $B_{j'} = c_{j' - 1}c_{j'}$ for all $j'$ by Lemma \ref{lem c}. Moreover, Lemma \ref{lem f} yields that $B_{j} = c_{j - 1} \vee K_{n_{j}} \vee c_{i}$. Let $c_{0} \in V(B_{1}) \setminus \{c_{1}\}$. We will show that $B_{1} = c_{0}c_{1}$. We let $a$ be a vertex in $V(B_{2}) - \{c_{1}, c_{2}\}$ if $j = 2$. Then, we let

\[ x        = \left\{\begin{array}{ll}
                     a       & if~~j = 2~~and\\
                     c_{2}   & if~~j > 2.
               \end{array}\right. \]
Consider $G + c_{0}x$. Since $c_{2}$ is a cut vertex of $G + c_{0}x$, $c_{2} \in D_{c_{0}x}$ by Lemma \ref{lem 2}(3). That is $x \in D_{c_{0}x}$ when $j > 2$. When $j = 2$, by Lemma \ref{lem f}, $xc_{2} \in E(G)$. Since $c_{2}\in D_{c_{0}x}$, by Lemma \ref{lem 1}(3), $x \in D_{c_{0}x}$. In both cases, $x \in D_{c_{0}x}$. If $|D_{c_{0}x} \cap (\cup^{k - 2}_{i = 2}V(B_{i}))| \leq k - 2$, then $(D_{c_{0}x} \cap (\cup^{k - 2}_{i = 2}V(B_{i}))) \cup \{c_{1}\} \succ_{c} G$ contradicting $\gamma_{c}(G) = k$. Therefore $|D_{c_{0}x} \cap (\cup^{k - 2}_{i = 2}V(B_{i}))| = k - 1$ by Lemma \ref{lem 1}(1). Thus $c_{1}, c_{0} \notin D_{c_{0}x}$ implying that $B_{1} = c_{0}c_{1}$. So $G \in \mathcal{G}_{1}(0, 0, ..., i_{j} = 1, ..., 0)$. This completes the proof.
\qed

\vskip 8 pt

\begin{Theorem}\label{thm g2}
For $k \geq 4$, $\mathcal{H}(0, 0, ..., 0, 3) = \mathcal{G}_{2}(k)$.
\end{Theorem}

\proof
By Lemma \ref{lem d}, $\mathcal{G}_{2}(k) \subseteq \mathcal{H}(0, 0, ..., 3)$. Thus, it is sufficient to show that $\mathcal{H}(0, 0, ..., 3) \subseteq \mathcal{G}_{2}(k)$.
\vskip 5 pt

\indent As $b_{i} = 0$ for all $2 \leq i \leq k - 3$, by Lemma \ref{lem c}, $B_{i} = c_{i - 1}c_{i}$. By Lemma \ref{lem pmp0} and similar arguments in Theorem \ref{thm giii}, we have that $B_{1} = c_{0}c_{1}$.
\vskip 5 pt

\vskip 5 pt

\indent We will show that $B_{k - 2}$ satisfies $(1)$ of the class $\mathcal{B}_{3}$. Recall that $\check{B}_{k - 2} = B_{k - 2} - c_{k - 3}$. Let $D'$ be a $\gamma_{c}$-set of $B_{k - 2}$. Suppose that $|D'| \leq 2$. To dominate $c_{k - 3}$, we have $D' \cap A \neq \emptyset$. Thus, $\{c_{1}, ..., c_{k - 3}\} \cup D' \succ_{c} G$. But $|\{c_{1}, ..., c_{k - 3}\} \cup D'| = k - 1$ contradicting the minimality of $k$. Therefore, $|D| = 3$. Thus, to prove that $B_{k - 2}$ satisfies $(1)$, it suffices to give a connected dominating set of size $3$ of $B_{k - 2}$ containing a chosen vertex from $\check{B}_{k - 2}$. For a vertex $v$ of $\check{B}_{k - 2}$, consider $G + c_{0}v$. Lemma \ref{lem 1}(2) implies that $\{c_{0}, v\} \cap D_{c_{0}v} \neq \emptyset$. Lemma \ref{lem 1}(1) implies also that $|D_{c_{0}v}| \leq k - 1$. We first show that $\{c_{0}\} \neq D_{c_{0}v} \cap \{c_{0}, v\}$. Suppose to the contrary that $\{c_{0}\} = D_{c_{0}v} \cap \{c_{0}, v\}$. Since $(G + uv)[D_{uv}]$ is connected, there exists $w \in V(\check{B}_{k - 2})$ which is adjacent to $v$. Because $D_{c_{0}v} \succ_{c} G + c_{0}v$, $w$ is adjacent to a vertex of $D_{uv} \cap V(B_{k - 2} - v)$. So
\vskip 5 pt

\indent $(D_{uv} - \{c_{0}\}) \cup \{w\} \succ_{c} G$.
\vskip 5 pt

\noindent This contradicts the minimality of $k$. Thus, $\{c_{0}\} \neq D_{c_{0}v} \cap \{c_{0}, v\}$. Therefore $\{c_{0}, v\} \subseteq D_{c_{0}v}$ or $\{v\} = D_{c_{0}v} \cap \{c_{0}, v\}$.
\vskip 5 pt

\noindent $\textbf{Case 1 :}$ $\{c_{0}, v\} \subseteq D_{c_{0}v}$\\
\noindent Let
\vskip 5 pt

\indent $i = max\{1 \leq j \leq k - 3 : G[\{c_{0}, c_{1}, c_{2}, ..., c_{j}\} \cap D_{c_{0}y}]$ is connected $\}$.
\vskip 5 pt

\noindent We first consider the case when $i = k - 3$. Thus $\{c_{1}, c_{2}, ..., c_{k - 3}\} \subseteq D_{c_{0}v}$. As $|D_{c_{0}v}| \leq k - 1$ and $\{c_{0}, v\} \subseteq D_{c_{0}v}$, we must have
\vskip 5 pt

\indent $D_{c_{0}v} = \{c_{0}, c_{1}, ..., c_{k - 3}, v\}$.
\vskip 5 pt

\noindent So $v \succ V(B_{k - 2}) - A$ and $N_{A}(v) = \emptyset$, otherwise $\{c_{1}, ..., c_{k - 3}, w, v\} \succ_{c} G$ where $w \in N_{A}(v)$, contradicting the minimality of $k$. Let $u \in N_{B_{k - 2}}(v)$ such that $u$ is adjacent to a vertex $a$ in $A$. Thus $\{v, u, a\} \succ_{c} B_{k - 2}$ and so $B_{k - 2}$ satisfies $(1)$.
\vskip 5 pt

\indent We now consider the case when $i = k - 4$. Let $D'_{c_{0}v} = D_{c_{0}v} \cap V(\check{B}_{k - 2})$. Clearly, $v \in D'_{c_{0}v}$. Since $\{c_{0}, c_{1}, ..., c_{k - 4}\} \subseteq D_{c_{0}v}$ and $|D_{c_{0}v}| \leq k - 1$, it follows that $|D'_{c_{0}v}| \leq 2$. If $|D'_{c_{0}v}| = 1$, then $D'_{c_{0}v}= \{v\}$ implying that $v \succ \check{B}_{k - 2}$, in particular, $v \succ A$. Thus, $\{c_{1}, ..., c_{k - 3}, w, v\} \succ_{c} G$ where $w \in N_{A}(v)$ contradicting the minimality of $D$. Hence, we let $D'_{c_{0}v} = \{v, v'\}$. Since $D_{c_{0}v} \succ_{c} G + c_{0}v$, $D'_{c_{0}v} \succ_{c} \check{B}_{k - 2}$. Hence, for a vertex $a$ in $A$, $D'_{c_{0}v} \cup \{a\} \succ_{c} B_{k - 2}$. Therefore, $B_{k - 2}$ satisfies $(1)$.
\vskip 5 pt

\indent We now consider the case when $i = k - 5$. Thus $\{c_{0}, c_{1}, ..., c_{k - 5}\} \subseteq D_{c_{0}v}$. So $|D'_{c_{0}v}| \leq 3$ and, $D'_{c_{0}v} \cap A \neq \emptyset$ to dominate $c_{k - 3}$. So, $B_{k - 2}$ satisfies $(1)$.
\vskip 5 pt

\indent We finally consider the case when $i \leq k - 6$. To dominate $c_{i + 2}$, we have that $c_{i + 3} \in D_{c_{0}v}$. By the connectedness of $(G + c_{0}v)[D_{c_{0}v}]$, $\{c_{i + 3}, ..., c_{k - 3}\} \subseteq D_{c_{0}v}$. Thus, $\{c_{0}, c_{1}, ..., c_{i}\} \cup \{c_{i + 3}, ..., c_{k - 3}\} \cup D'_{c_{0}v} \subseteq D_{c_{0}v}$ implying that
 \begin{center}
 $k - 4 + |D'_{c_{0}v}| = (i + 1) + ((k - 3) - (i + 3) + 1) + |D'_{c_{0}v}| \leq k - 1$.
\end{center}
\noindent Therefore, $|D'_{c_{0}v}| \leq 3$. To dominate $c_{k - 3}$, $D'_{c_{0}v} \cap A \neq \emptyset$. Thus $B_{k - 2}$ satisfies $(1)$ and this completes the proof of Case 1.
\vskip 5 pt

\noindent $\textbf{Case 2 :}$ $\{v\} = D_{c_{0}v} \cap \{c_{0}, v\}$\\
\indent To dominate $c_{1}$, we have that $c_{2} \in D_{c_{0}v}$. By the connectedness of $(G + c_{0}v)[D_{c_{0}v}]$, $\{c_{2}, c_{3}, ...,$ $c_{k - 3}\} \subseteq D_{c_{0}v}$ and $D_{c_{0}v} \cap A \neq \emptyset$. As $|D_{c_{0}v}| \leq k - 1$, we must have $|D_{c_{0}v} - \{c_{2}, c_{3}, ..., c_{k - 3}\}| \leq 3$. So $B_{k - 2}$ satisfies $(1)$. This completes the proof of Case 2.
\vskip 5 pt

\indent We finally show that $B_{k - 2}$ satisfies $(2)$ of the class $\mathcal{B}_{2}$. Let $x$ and $y$ be non-adjacent vertices of $B_{k - 2}$. Lemma \ref{lem 1}(2) implies that $\{x, y\} \cap D_{xy} \neq \emptyset$. Lemma \ref{lem 1}(1) implies also that $|D_{xy}| \leq k - 1$. To dominate $c_{0}$, we have that $c_{1} \in D_{xy}$. Let $D^{B}_{xy} = D_{xy} \cap V(\check{B}_{k - 2})$. By the connectedness of $(G + xy)[D_{xy}]$, $D^{B}_{xy} \cap A \neq \emptyset$ and $\{c_{1}, c_{2}, ..., c_{k - 3}\} \subseteq D_{xy}$. As $|D_{xy}| \leq k - 1$, we must have $|D^{B}_{xy}| = |D_{xy} \cap V(\check{B}_{k - 2})| = |D_{xy} - \{c_{1}, c_{2}, ..., c_{k - 3}\}| \leq 2$. Hence, $B_{k - 2}$ satisfies $(2)$. Therefore $B_{k - 2} \in \mathcal{B}_{3}$. This completes the proof.
\qed

\vskip 5 pt

\indent No, we are ready to prove our first main result, Theorem \ref{thm k-3}. For completeness, we restate the theorem
\vskip 5 pt

\noindent \textbf{Theorem \ref{thm k-3}}. \emph{For an integer $k \geq 4$, $\mathcal{Z}(k, k - 3) = \mathcal{G}_{1}(i_{1}, i_{2}, ..., i_{k - 3}) \cup \mathcal{G}_{2}(k)$.}
\vskip 5 pt

\proof
In view of Lemma \ref{lem e2}, $\mathcal{Z}(k, k - 3) = \mathcal{H}(0, 0, 0, ..., 3) \cup \mathcal{H}(b_{1}, b_{2}, ..., b_{k - 3}, 2)$ where $b_{i} = 1$ for some $1 \leq i \leq k - 3$ and $b_{j} = 0$ for all $1 \leq j \neq i \leq k - 3$. Moreover, Theorems \ref{thm giii} and \ref{thm g2} imply that $\mathcal{Z}(k, k - 3) = \mathcal{G}_{1}(i_{1}, i_{2}, ..., i_{k - 3}) \cup \mathcal{G}_{2}(k)$. This completes the proof.
\qed

\subsection{$k$-$\gamma_{c}$-Critical Graphs with Prescribed Cut Vertices}\label{real}
\noindent In this section, we prove Theorem \ref{thm:main111}. First, we introduce structure of some subgraphs. For any block $H$ of a graph $G$, $H$ is call a \emph{block} $\mathcal{H}_{\ell}$ for $\ell \geq 2$ if $H$ consists of a vertex $x$ and $U_{1}, U_{2}, ..., U_{\ell}$ as vertex sets of order at least $2$ and
\begin{itemize}
  \item $G[\{x\} \cup U_{1}]$, $G[U_{1} \cup U_{2}], G[U_{2} \cup U_{3}], ..., G[U_{\ell - 2} \cup U_{\ell - 1}]$ and $G[U_{\ell}]$ are complete.
  \item For each $u \in V(U_{\ell - 1})$, $|N_{U_{\ell}}(u)| = |U_{\ell}| - 1$.
  \item For each $u' \in V(U_{\ell})$, there exists $u \in V(U_{\ell - 1})$ such that $uu' \in E(G)$.
\end{itemize}

\noindent We say that $x$ is the \emph{head} of a block $\mathcal{H}_{\ell}$.
\vskip 5 pt

\noindent Let\\
\indent $\mathcal{D}(k, \zeta, \zeta_{0}) :$ the class of all $k$-$\gamma_{c}$-critical graphs with $\zeta$ cut vertices containing a block $B$ such that $\zeta(B) = \zeta_{0}$.
\vskip 5 pt

\indent We next introduce the following class that we use to establish the existence of graphs in $\mathcal{D}(k, \zeta, \zeta_{0})$.
\vskip 5 pt

\noindent \textbf{The class $\mathcal{F}(p, q, r)$}:\\
\indent Let $H_{1}, H_{2}, ..., H_{p}$ be $p$ of $\mathcal{H}_{2}$ blocks. Further, we let $H_{p + 1}$ be $P_{q}$ (a path of $q$ vertices), and let $H_{p + 2}$ be an $\mathcal{H}_{r}$ block. Let $c_{i}$ be the head of $H_{i}$ for $1 \leq i \leq p + 2$. A graph $G \in \mathcal{F}(p, q, r)$ is obtained from $H_{1}, ..., H_{p + 2}$ by joining edges between vertices in $\{c_{i} : 1 \leq  i \leq p + 2\}$ to form a clique.
\vskip 5 pt

\begin{Lemma}\label{lem pqr}
If a graph $G$ is in the class $\mathcal{F}(p, q, r)$, then $G$ is $(r + q + 3p)$-$\gamma_{c}$-critical with $p + q$ cut vertices having a block that contains $p + 2$ cut vertices.
\end{Lemma}

\proof
For $i \in [p]$, let $U^{i}_{1} = N_{H_{i}}(c_{i})$ and $U^{i}_{2} = V(H_{i}) - U^{i}_{1}$. We see that $U^{i}_{1}$ and $U^{i}_{2}$ are, respectively, the same as $U_{1}$ and $U_{2}$ of the block $\mathcal{H}_{2}$. Let $a_{i} \in U^{i}_{1}$. By the construction of the block $\mathcal{H}_{2}$, there exists $b_{i} \in U^{i}_{2}$ such that $a_{i}b_{i} \in E(G)$ and $a_{i}$ is not adjacent to the vertex $a'_{i} \in U^{i}_{2}$. Let $H_{p + 1} = d_{0}, d_{1}, ..., d_{q - 1}$ where $c_{p + 1} = d_{0}$. For the block $H_{p + 2}$, let $\tilde{U}_{i}$ be a vertex subset of $V(H_{p + 2})$ which is the same as $U_{i}$ of the block $\mathcal{H}_{r}$ and $u_{i} \in \tilde{U}_{i}$ for $i \in [r]$.
\vskip 5 pt

\indent Lemma \ref{lem 2}(3) implies that $c_{i} \in D$ for $i \in [p + 2]$. Consider $H_{p + 2}$. By the connectedness of $G[D]$, $|\tilde{U}_{i} \cap D_{c}| \geq 1$ for $i \in [r - 1]$. Because every vertex in $\tilde{U}_{r - 1}$ does not dominate $\tilde{U}_{r}$, $|D \cap (\tilde{U}_{r - 1} \cup \tilde{U}_{r})| \geq 2$. Since $c_{p + 2} \in D$, $|D \cap V(H_{p + 2})| \geq r + 1$. By the same arguments, $|D \cap V(H_{i})| \geq 3$ for $i = 1, 2, ..., p$. Finally, to dominate $d_{q - 1}$, $\{c_{p + 1}, d_{1}, d_{2}, ..., d_{q - 2}\} \subseteq D$. So $|D| \geq (r + 1) + 3p + (q - 1) = r + q + 3p$. It is not difficult to see that $(\cup^{p}_{i = 1}\{a_{i}, b_{i}\}) \cup \{c_{1}, c_{2}, ..., c_{p + 2}\} \cup \{d_{1}, .., d_{q - 2}\} \cup \{u_{1}, u_{2}, ..., u_{r}\} \succ_{c} G$. So $|D| \leq 2p + (p + 2) + (q - 2) + r = r + q + 3p$. Thus $\gamma_{c}(G) = r + q + 3p$.
\vskip 5 pt

\indent Let $k = r + q + 3p$. For a pair of non adjacent vertices $u, v\in V(G)$. Consider $G + uv$. To establish the criticality, it suffices to show that there exists a dominating set $D_{uv}$ of $G + uv$ containing less than $k$ vertices. For $i \in [p + 2]$, let $D_{i} = D \cap V(H_{i})$. We distinguish $3$ cases.
\vskip 5 pt

\noindent \textbf{Case 1 :} $\{u, v\} \cap (\cup^{p}_{i = 1}H_{i}) \neq \emptyset$.\\
\indent Without loss of generality, let $u \in V(H_{1})$.
\vskip 5 pt

\noindent \textbf{Subcase 1.1 :} $u = c_{1}$. Clearly, $v \notin \{c_{2}, c_{3}, ..., c_{p + 2}\}$.
\begin{itemize}
  \item If $v \in U^{1}_{2}$, then let $D_{uv} = \{u, v\} \cup (\cup^{p + 2}_{i = 2}D_{i})$.
\end{itemize}
\noindent For $j \geq 2$,
\begin{itemize}
  \item if $v \in U^{j}_{1}$, then let $D_{uv} = (\cup_{i \neq j}D_{i}) \cup \{v, w\}$ where $w \in N_{U^{j}_{2}}(v)$, and
  \item if $v \in U^{j}_{2}$, then let $D_{uv} = (\cup_{i \neq j}D_{i}) \cup \{v, c_{j}\}$.
\end{itemize}

\noindent We now consider the case when $v \in V(H_{p + 1})$. We let $D_{uv} = (\cup_{i \neq p + 1}D_{i}) \cup \{d_{1}, d_{2}, ..., $ $v, ..., d_{q - 2}\}$ if $v \neq d_{q - 1}$ and let $D_{uv} = D - \{d_{q - 2}\}$ if $v = d_{q - 1}$. When $v \in V(H_{p + 2})$, we let $D_{uv} = (\cup_{i \neq p + 2}D_{i}) \cup \{u_{1}, u_{2}, ..., v, ..., u_{r}\}$ where $u_{1}, u_{2}, ..., v, ..., u_{r}$ is a path such that $v = u_{j}$ for some $j \in [r]$.
\vskip 5 pt

\noindent \textbf{Subcase 1.2 :} $u \in U^{1}_{1}$. Without loss of generality let $u = a_{1}$. So $ua'_{1} \notin E(G)$ and $ub_{1} \in E(G)$. Clearly, $v \neq c_{1}$. If $v = c_{i}$ for some $\{2, ..., p\}$, then we can find $D_{uv}$ by the same arguments as Subcase 1.1. Thus, we may consider when $v \neq c_{i}$ for $i = 1, 2, ..., p$. If $v = a'_{1}$, then let $D_{uv} = (\cup_{i \neq 1}D_{i}) \cup \{a_{1}, c_{1}\}$.
\noindent For $j \geq 2$,
\begin{itemize}
  \item if $v \in U^{j}_{1}$, then let $D_{uv} = (\cup_{i \neq 1, j}D_{i}) \cup \{c_{1}, u, b_{1}, v, w\}$ where $w \in N_{U^{j}_{2}}(v)$, and
  \item if $v \in U^{j}_{2}$, then let $D_{uv} = (\cup_{i \neq 1, j}D_{i}) \cup \{c_{1}, u, b_{1}, v, c_{j}\}$.
\end{itemize}

\noindent We now consider the case when $v \in V(H_{p + 1}) \cup V(H_{p + 2})$. If $v \in \{c_{p + 1}, c_{p + 2}\}$, then $D_{uv} = (\cup_{i \neq 1}D_{i}) \cup \{u, b_{1}\}$. Thus, we may assume that $v \notin \{c_{p + 1}, c_{p + 2}\}$. We let $\tilde{D} = (\cup^{p}_{i = 2}D_{i}) \cup \{c_{1}, u, b_{1}\}$. If $v = d_{i}$ for some $i \in [q - 2]$, then we let $D_{uv} = \tilde{D} \cup \{d_{1}, ..., d_{q - 2}\} \cup D_{p + 2}$. If $v = d_{q - 1}$, then we let $D_{uv} = \tilde{D} \cup \{d_{0}, ..., d_{q - 3}\} \cup D_{p + 2}$. Finally, if $v = u_{i} \in \tilde{U}_{i}$ for some $i \in [r]$, then we let $D_{uv} = \tilde{D} \cup D_{p + 1} \cup \{u_{1}, u_{2}, ..., v, ..., u_{r}\}$.
\vskip 5 pt

\noindent \textbf{Subcase 1.3 :} $u \in U^{1}_{2}$. Without loss of generality let $u = b_{1}$. Clearly $ua_{1} \in E(G)$. By the same arguments as Subcases 1.1 and 1.2, we consider only when $v \notin \cup^{p}_{i = 2}(\{c_{i}\} \cup U^{i}_{1})$. If $v \in U^{1}_{1}$, then let $D_{uv} = \{v, c_{1}\} \cup (\cup_{i \neq 1}D_{i})$. If $v \in U^{j}_{2}$, then let $D_{uv} = (\cup_{i \neq 1, j}D_{i})$ $\cup \{u, v, a_{1}, c_{1}\}$. If $v \in \{c_{p + 1}, c_{p + 2}\}$, then $D_{uv} = (\cup_{i \neq 1}D_{i}) \cup \{u, c_{1}\}$. Thus, we may assume that $v \in V(H_{p + 1}) \cup V(H_{p + 2}) - \{c_{p + 1}, c_{p + 2}\}$. We let $\tilde{D} = D_{uv} = (\cup^{p}_{i = 2}D_{i}) \cup \{c_{1}, u, a_{1}\}$ and we can find $D_{uv}$ by same arguments as Subcase 1.2.
\vskip 5 pt

\noindent \textbf{Case 2 :} $\{u, v\} \cap (\cup^{p}_{i = 1}V(H_{i})) = \emptyset$ and $\{u, v\} \cap V(H_{p + 1}) \neq \emptyset$.\\
\indent Without loss of generality let $u \in V(H_{p + 1})$.
\vskip 5 pt

\noindent \textbf{Subcase 2.1 :} $u = c_{p + 1}$. Suppose first that $v \in V(H_{p + 1})$. Thus $v = d_{j}$ for some $j > 2$. If $j = 3$, then let $D_{uv} = (\cup_{i \neq p + 1}D_{i}) \cup \{c_{p + 1}, d_{3}, d_{4}, ..., d_{q - 2}\}$. If $3 < j < q - 1$, then $D_{uv} = (\cup_{i \neq p + 1}D_{i}) \cup \{c_{p + 1}, d_{4},..., v ..., d_{q - 2}\}$. If $j = q - 1$, then let $D_{uv} = D - \{d_{q - 2}\}$. If $v \in V(H_{p + 2})$, then we let $D_{uv} = (\cup_{i \neq p + 2}D_{i}) \cup \{u_{1}, u_{2}, ..., v, ..., u_{r}\}$.
\vskip 5 pt

\noindent \textbf{Subcase 2.2 :} $u \in V(H_{p + 1}) - \{c_{p + 1}\}$. Suppose that $v \in V(H_{p + 1})$. By the same arguments as Subcase 2.1, $v \in V(H_{p + 1}) - \{c_{p + 1}\}$. Without loss of generality let $u = d_{j}$ and $v = d_{j'}$ such that $j < j'$. If $j' < q - 1$, then let $D_{uv} = (\cup_{i \neq p + 1}D_{i}) \cup \{c_{p + 1}, d_{1}, d_{2}, ..., d_{j - 1}, u, d_{j + 2}, ..., v, ..., d_{q - 2}\}$. If $j' = q - 1$, then let $D_{uv} = (\cup_{i \neq p + 1}D_{i})$ $\cup \{c_{p + 1}, d_{1}, ..., u, ..., d_{q - 3}\}$. Suppose that $v \in V(H_{p + 2})$. If $v = c_{p + 2}$ and $j < q - 1$, then let $D_{uv} = (\cup_{i \neq p + 1}D_{i}) \cup \{d_{1}, ..., u, ..., d_{q - 2}\}$. If $v = c_{p + 2}$ and $j = q - 1$, then let $D_{uv} = D - \{d_{q - 2}\}$. We now consider the case when $v \in V(H_{p + 2}) - \{c_{p + 2}\}$. Let $\tilde{D}_{p + 2} = \{c_{p + 2}, u_{1}, ..., u_{r}\}$ and $\tilde{D}_{p + 1} = \{d_{1}, d_{2}, ..., d_{q - 2}\}$. Let $D_{uv} = (\cup^{p}_{i = 1}D_{i}) \cup (\tilde{D}_{p + 1} \cup \tilde{D}_{p + 2})$ by the same arguments as Subcase 2.1.
\vskip 5 pt

\noindent \textbf{Case 3 :} $\{u, v\} \subseteq V(H_{p + 2})$. Without loss of generality let $u \in \tilde{U}_{j}$ and $v \in \tilde{U}_{j'}$ where $j < j'$. Let $c_{p + 2} = U_{0}$. If $u \notin \tilde{U}_{r - 1}$, then there exists a path $u_{0}, u_{1}, ..., u_{r}$ such that $u = u_{j}$ and $v = u_{j'}$ where $u_{i} \in \tilde{U}_{i}$ for $0 \leq i \leq r$. Let $D_{uv} = (\cup_{i \neq p + 2}D_{i}) \cup \{u_{0}, u_{1}, ..., u_{j - 1}, u, u_{j + 2}, ..., v, ..., u_{r}\}$. If $u \in \tilde{U}_{r - 1}$, then $v$ is the only one vertex in $\tilde{U}_{r}$ which $u$ is not adjacent. Thus $D_{uv} = (\cup_{i \neq p + 2}D_{i}) \cup \{U_{0}, U_{1}, ..., U_{r - 2}, u\}$.
\vskip 5 pt

\indent Finally, we see that $c_{1}, ..., c_{p + 2}, d_{1}, ..., d_{q - 2}$ are all the cut vertices of $G$. Thus, $G$ has $p + q$ cut vertices. Further, the block $G[\{c_{1}, ..., c_{p + 2}\}]$ has $c_{1}, ..., c_{p + 2}$ as the cut vertices of $G$. Therefore, there is a block containing $p + 2$ cut vertices. This completes the proof.
\qed
\vskip 5 pt

\noindent Now, we are ready to prove Theorem \ref{thm:main111}. For completeness, we restate the theorem.
\vskip 5 pt

\noindent \textbf{Theorem \ref{thm:main111}}.
\emph{For all $k \geq 4$, $2 \leq \zeta \leq k - 2$ and $2 \leq \zeta_{0} \le \min \left\{ \left\lfloor \frac{k + 2}{3} \right\rfloor, \zeta \right\}$, there exists a $k$-$\gamma_{c}$-critical graph with $\zeta$ cut vertices having a block that contains $\zeta_{0}$ cut vertices, namely, $\mathcal{D}(k, \zeta, \zeta_{0}) \neq \emptyset$.}

\proof
In view of Lemma \ref{lem pqr}, for all $k \geq 4$, $2 \leq \zeta \leq k - 2$ and $2 \leq \zeta_{0} \le \min \left\{ \left\lfloor \frac{k + 2}{3} \right\rfloor, \zeta \right\}$, we have $\mathcal{F}(\zeta_{0} - 2, \zeta - \zeta_{0} + 2, k - \zeta - 2\zeta_{0} + 4) \subseteq \mathcal{D}(k, \zeta, \zeta_{0})$. This completes the proof.
\qed

\subsection{Factor Criticality of $k$-$\gamma_{c}$-Critical Graphs}\label{fact}
\noindent In this section, for $k \geq 3$, we will use the property of graphs in the class $\mathcal{P}(k)$ which is given in Section \ref{pre}. First, we may prove that the class $\mathcal{P}(k)$ is non-empty for $k \geq 3$.




\vskip 5 pt

\begin{Lemma}\label{lem 5}
For all $k \geq 3$, $\mathcal{P}(k) \neq \emptyset$.
\end{Lemma}

\proof
For an integer $k \geq 3$, we let $C_{k + 2} = c_{1}, c_{2}, ..., c_{k + 2}, c_{1}$ be a cycle of length $k + 2$. It is well known that $C_{k + 2}$ is a $k$-$\gamma_{c}$-critical graph. In this proof, all subscripts are taken modulo $k + 2$. Observe that $C_{k + 2}[\{c_{j}, c_{j + 1}\}]$ is a maximal complete subgraph for all $1 \leq j \leq k + 2$. Thus, without loss of generality, it suffices to show that complete subgraph $C_{k + 2}[\{c_{1}, c_{2}\}]$ satisfies $(i)$ and $(ii)$ of $\mathcal{P}(k)$.
\vskip 5 pt

\indent For a vertex $c_{j}$ where $1 \leq j \leq k + 2$, let

\begin{center}
$S_{j} = \{c_{j}, c_{j + 1}, ..., c_{k + 2}, c_{1}, ..., c_{j - 3}\} \succ_{c} C_{k + 2}$.
\end{center}

\noindent Clearly, $S_{j} \cap \{c_{1}, c_{2}\} \neq \emptyset$ if $j \neq 3$. Thus we consider the case when $j = 3$. In this case $c_{k + 2} \in S_{3}$. Let $S' = (S_{3} - \{c_{k + 2}\}) \cup \{c_{2}\}$. Hence, $S' \succ_{c} C_{k + 2}$. Moreover, $S' \cap \{c_{1}, c_{2}\} \neq \emptyset$. Therefore $C_{k + 2}[\{c_{1}, c_{2}\}]$ satisfies $(i)$.
\vskip 5 pt

\indent We now let $c_{j}$ and $c_{l}$ be non-adjacent vertices of $C_{k + 2}$. So $|j - l| \geq 2$. Consider $C_{k + 2} + c_{j}c_{l}$. We partition $V(C_{k + 2})$ to
\begin{center}
$C^{1} = \{c_{j}, c_{j + 1}, ..., c_{l - 1}\}$ and $C^{2} = \{c_{l}, c_{l + 1}, ..., c_{j - 1}\}$.
\end{center}
As $k + 2 \geq 5$, at least one of $C^{1}$ or $C^{2}$ must have at least three vertices. Without loss of generality let it be $C^{1}$. We further let
\vskip 5 pt
\indent
\[ D  = \left\{\begin{array}{ll}
                     \{c_{j}, c_{j + 1}, ..., c_{l - 3}\} \cup \{c_{l}, c_{l + 1}, ..., c_{j - 3}\}, & if~~|C^{2}| \geq 3~~and\\
                     \{c_{j}, c_{j + 1}, ..., c_{l - 2}\}, & otherwise.
               \end{array}\right. \]
Clearly $D \succ_{c} C_{k + 2}$ and $|D| < k$. We first consider the case when $|C^{2}| \geq 3$. Note that $D$ in this case contains $k - 2$ vertices. If $\{c_{1}, c_{2}\}$ is neither $\{c_{l - 1}, c_{l - 2}\}$ nor $\{c_{j - 1}, c_{j - 2}\}$, then $D \cap \{c_{1}, c_{2}\} \neq \emptyset$. This implies that $C_{k + 2}[\{c_{1}, c_{2}\}]$ satisfies $(ii)$. Thus we may assume that $\{c_{1}, c_{2}\} = \{c_{l - 1}, c_{l - 2}\}$. As $|D| = k - 2$, we must have $|D \cup \{c_{l - 1}\}| < k$. Moreover $D \cup \{c_{l - 1}\} \succ_{c} C_{k + 2} + c_{j}c_{l}$ because $D \succ_{c} + c_{j}c_{l}$. Hence $C_{k + 2}[\{c_{1}, c_{2}\}]$ satisfies $(ii)$. The case $\{c_{1}, c_{2}\} = \{c_{j - 1}, c_{j - 2}\}$ can be proved by the same arguments.
\vskip 5 pt

\indent We now consider the case when $|C^{2}| < 3$, in fact $C^{2} = \{c_{l}, c_{l + 1}\}$ and $c_{l + 2} = c_{j}$. If $\{c_{1}, c_{2}\}$ is neither $\{c_{l - 1}, c_{l}\}$ nor $\{c_{l}, c_{l + 1}\}$, then $D \cap \{c_{1}, c_{2}\}$ and this implies that $C_{k + 2}[\{c_{1}, c_{2}\}]$ satisfies $(ii)$. Thus we may assume that $\{c_{1}, c_{2}\}$ is $\{c_{l - 1}, c_{l}\}$ or $\{c_{l}, c_{l + 1}\}$. We observe that $c_{1} = c_{l - 1}$ and $c_{2} = c_{l}$ when $\{c_{1}, c_{2}\} = \{c_{l - 1}, c_{l}\}$, moreover, $c_{1} = c_{l}$ and $c_{2} = c_{l + 1}$ when $\{c_{1}, c_{2}\} = \{c_{l}, c_{l + 1}\}$. In both cases, $c_{l} \in \{c_{1}, c_{2}\}$. Let $D' = (D - \{c_{l - 2}\}) \cup \{c_{l}\}$. We see that $D' \succ_{c} C_{k + 2} + c_{j}c_{l}$. So $D' \cap \{c_{1}, c_{2}\} \neq \emptyset$ because $c_{l} \in \{c_{1}, c_{2}\}$ and $c_{l} \in D'$. Therefore $C_{k + 2}[\{c_{1}, c_{2}\}]$ satisfies $(ii)$. Hence $C_{k + 2} \in \mathcal{P}(k)$ and this completes the proof.
\qed

\vskip 5 pt

\indent In the following, we show how to apply the construction of some graphs in the class $\mathcal{P}(k)$ to establish the existence of $(k + l)$-$\gamma_{c}$-critical graphs with some property. For an integer $k \geq 1$ and $\ell \geq 1$, we let
\vskip 5 pt

\indent $\mathcal{Q}(k, \ell) :$ the class of $k$-$\gamma_{c}$-critical graphs $G$  with $\delta \geq \ell + 1$ such that\\
\indent $\quad \quad \quad \quad$ $G$ is not $\ell$-factor critical.
\vskip 5 pt

\noindent Hence, we may rewrite Observation \ref{obs 4c} in term of the class $\mathcal{Q}(k, 1)$.
\vskip 5 pt

\begin{ob}\label{obs 5}
$\mathcal{Q}(k, 1) = \emptyset$ for $k \in [2]$.
\end{ob}
\vskip 5 pt

\noindent For $k = 3$, Figure \ref{f:fig001} shows that there exists a $3$-$\gamma_{c}$-critical graph of odd order and $\delta \geq 2$ which is non-factor critical. Thus $\mathcal{Q}(3, 1) \neq \emptyset$. In the following, for $k \geq 4$, we show further that there exists a $k$-$\gamma_{c}$-critical graph which is non-factor critical.
\vskip 10 pt

\noindent \textbf{The class} $\mathcal{X}(s)$\\
\indent For an integer $s \geq 3$, let $A = \{a_{1}, a_{2}, ..., a_{s}\}$ and $B = \{b_{1}, b_{2}, ..., b_{s}\}$ be two disjoint sets of vertices. We further let $K_{s}$ be a copy of a complete graph of order $s$ such that $V(K_{s}) = \{y_{1}, y_{2},..., y_{s}\}$. A graph $G$ in the class $\mathcal{X}(s)$ can be constructed from $A, B$ and $K_{s}$ by adding edges according to the join operations that, for $1 \leq i \leq s$,
\begin{itemize}
  \item $a_{i} \vee (B - \{b_{i}\})$ and
  \item $a_{i} \vee (K_{s} - y_{i})$.
\end{itemize}
\noindent A graph in this class is illustrated by Figure \ref{X}.

\vskip 20 pt

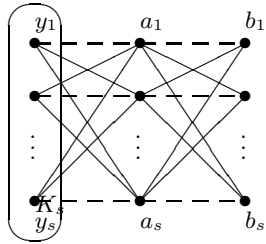
\begin{figure}[htb]
\begin{center}
\setlength{\unitlength}{0.7cm}
\begin{picture}(12, 6)
\put(3, 5){\circle*{0.2}}
\put(3, 4){\circle*{0.2}}
\put(3, 2){\circle*{0.2}}
\put(3, 3.5){\oval(1, 4.5)}
\put(5, 5){\circle*{0.2}}
\put(5, 4){\circle*{0.2}}
\put(5, 2){\circle*{0.2}}
\put(7, 5){\circle*{0.2}}
\put(7, 4){\circle*{0.2}}
\put(7, 2){\circle*{0.2}}

\put(3, 2){\line(2, 3){2}}
\put(3, 2){\line(1, 1){2}}
\put(5, 2){\line(2, 3){2}}
\put(5, 2){\line(1, 1){2}}
\put(3, 4){\line(1, -1){2}}
\put(3, 4){\line(2, 1){2}}
\put(5, 4){\line(1, -1){2}}
\put(5, 4){\line(2, 1){2}}
\put(3, 5){\line(2, -1){2}}
\put(3, 5){\line(2, -3){2}}
\put(5, 5){\line(2, -1){2}}
\put(5, 5){\line(2, -3){2}}

\multiput(3, 2)(0.5, 0){8}{\line(1,0){0.29}}
\multiput(3, 4)(0.5, 0){8}{\line(1,0){0.29}}
\multiput(3, 5)(0.5, 0){8}{\line(1,0){0.29}}

\put(6.9, 2.8){\footnotesize$\vdots$}
\put(4.9, 2.8){\footnotesize$\vdots$}
\put(2.9, 2.8){\footnotesize$\vdots$}
\put(5, 5.3){\footnotesize$a_{1}$}
\put(5, 1.5){\footnotesize$a_{s}$}
\put(7, 5.3){\footnotesize$b_{1}$}
\put(7, 1.5){\footnotesize$b_{s}$}
\put(3, 5.3){\footnotesize$y_{1}$}
\put(3, 1.5){\footnotesize$y_{s}$}
\put(3, 1.8){\footnotesize$K_{s}$}
\end{picture}
\caption{A graph $G$ in the class $\mathcal{X}(s)$}
\label{X}
\end{center}
\end{figure}
\vskip 10 pt

\indent The following lemma gives that $\mathcal{X}(s) \subseteq \mathcal{P}(4) \cap \mathcal{Q}(4, 1)$ for integer $s \geq 3$.
\vskip 5 pt

\begin{Lemma}\label{lem 6}
For an integer $s \geq 3$, $\mathcal{X}(s) \subseteq \mathcal{P}(4) \cap \mathcal{Q}(4, 1)$.
\end{Lemma}

\proof
For a given $s \geq 3$, let $G$ be in the class $\mathcal{X}(s)$. We first show that $\gamma_{c}(G) = 4$. Suppose to the contrary that there exists a connected dominating set $D$ of size less than $4$. We first consider the case when $D \cap V(K_{s}) \neq \emptyset$. To dominate $B$, $|D \cap A| \geq 2$. Therefore $D = \{y_{i}, a_{j}, a_{l}\}$. By the connectedness of $G[D]$, $i \notin \{j, l\}$. Hence $D$ does not dominate $a_{i}$, a contradiction. Thus, we consider the case when $D \cap V(K_{s}) = \emptyset$. So, to dominate $K_{s}$, $|D \cap A| \geq 2$. Since $A$ is an independent set, by the connectedness of $G[D]$, $|D \cap B| \geq 1$. Therefore $D = \{a_{j}, a_{l}, b_{i}\}$. Similarly $i \notin \{j, l\}$ and this implies that $D$ does not dominate $a_{i}$, a contradiction. Thus $\gamma_{c}(G) \geq 4$.
\vskip 5 pt

\indent We observe that $K_{s}$ is a maximal complete subgraph of $G$. We do not only show that $\gamma_{c}(G) \leq 4$ but we also show that, for a vertex $a$ of $G$, there exists a connected dominating set $D_{a}$ of $G$ containing $a$ and $D \cap V(K_{s}) \neq \emptyset$. That is we show that $K_{s}$ satisfies $(i)$ of graphs in the class $\cP(k)$. For $1 \leq i \neq i' \neq i'' \leq s$, we have
\vskip 5 pt

\indent $D_{y_{i}} = D_{a_{i}} = \{y_{i}, y_{i'}, a_{i}, a_{i'}\}$ and
\vskip 5 pt

\indent $D_{b_{i}} = \{b_{i}, a_{i'}, y_{i''}, a_{i}\}$
\vskip 5 pt

\noindent Hence $\gamma_{c}(G) \leq 4$ and thus $\gamma_{c}(G) = 4$. Moreover, $K_{s}$ satisfies $(i)$.
\vskip 5 pt

\indent We finally establish the criticality. Further, we show that, for non-adjacent vertices $x$ and $y$ of $G$, $D_{xy} \cap V(K_{s}) \neq \emptyset$. That is we will show that $K_{s}$ satisfies $(ii)$ of graphs in the class $\cP(k)$. We first consider the case when $\{x, y\} \cap B \neq \emptyset$. If $\{x, y\} = \{b_{i}, b_{j}\}$, then $D_{xy} = \{b_{i}, a_{j}, y_{l}\}$. If $\{x, y\} = \{a_{i}, b_{i}\}$, then $D_{xy} = \{a_{i}, y_{j}, y_{l}\}$. If $\{x, y\} = \{b_{i}, y_{j}\}$ (in this case $y_{j}$ could be $y_{i}$), then $D_{xy} = \{y_{j}, y_{l}, a_{i}\}$. We now consider the case when $\{x, y\} \cap B = \emptyset$. Thus $\{x, y\}$ is either $\{a_{i}, a_{j}\}$ or $\{a_{i}, y_{i}\}$. In both cases, $D_{xy} = \{a_{i}, a_{j}, y_{i}\}$. Thus $G$ is a $4$-$\gamma_{c}$-critical graph, in particular, $G \in \mathcal{P}(4)$.

\noindent Finally, let $s$ be odd number and $S = A$. Thus, $\omega_{o}(G - S)$ has $K_{s}$ and $b_{1}, ..., b_{s}$ as $s + 1$ odd components. Thus $\omega_{o}(G - S) = s + 1 > s - 1 = |S| - 1$. By Theorem \ref{thm favaron}, $G$ is non-factor critical. Thus, $G \in \mathcal{Q}(4, 1)$ and this completes the proof.
\qed

\vskip 5 pt

\indent We will use a graph in the class $\mathcal{P}(4)$ to show that $\mathcal{Q}(k, 1) \neq \emptyset$ for all $k \geq 6$. For $k = 5$, we also provide a graph $G_{5}(l_{1}, l_{2})$ in the class $\mathcal{Q}(5, 1)$ by the following construction. Let $u, x, y, z$ and $w$ be five different vertices. We also let $P_{2} = x', y'$ be a path of length one and $K_{l_{1}}, K_{l_{2}}$ be two copies of complete graphs of order $l_{1} \geq 2$ and $l_{2} \geq 2$ respectively, moreover, $l_{1} + l_{2}$ is even number. The graph $G_{5}(l_{1}, l_{2})$ is constructed by adding edges according to the join operations :
\begin{itemize}
  \item $u \vee K_{l_{1}} \vee K_{l_{2}} \vee P_{2}$
  \item $x \vee x', y \vee y'$
  \item $w \vee P_{2}$ and
  \item $z \vee \{x, y, w\}$.
\end{itemize}
\noindent Figure \ref{G5} illustrates the graph $G_{5}(l_{1}, l_{2})$.
\vskip 5 pt

\vskip 20 pt
\begin{figure}[htb]
\begin{center}
\setlength{\unitlength}{0.7cm}
\begin{picture}(20, 5)
\put(3.5, 3){\circle*{0.2}}
\put(3.5, 3){\line(1, 1){1}}
\put(3.5, 3){\line(1, -1){1}}
\put(5, 3){\oval(1, 3)}
\put(5.6, 2.5){\line(1, 0){0.7}}
\put(5.6, 3.5){\line(1, 0){0.7}}
\put(5.6, 3){\line(1, 0){0.7}}
\put(7, 3){\oval(1, 3)}

\put(8.8, 2){\circle*{0.2}}
\put(8.8, 4){\circle*{0.2}}

\put(8.8, 2){\line(-1, 0){1.3}}
\put(8.8, 2){\line(-2, 3){1.3}}
\put(8.8, 4){\line(-1, 0){1.3}}
\put(8.8, 4){\line(-2, -3){1.3}}
\put(8.8, 2){\line(0, 1){2}}

\put(10.8, 1){\circle*{0.2}}
\put(10.8, 3){\circle*{0.2}}
\put(10.8, 5){\circle*{0.2}}
\put(8.8, 2){\line(2, 1){2}}
\put(8.8, 2){\line(2, -1){2}}
\put(8.8, 4){\line(2, 1){2}}
\put(8.8, 4){\line(2, -1){2}}

\put(12.8, 3){\circle*{0.2}}
\put(12.8, 3){\line(-1, 0){2}}
\put(12.8, 3){\line(-1, 1){2}}
\put(12.8, 3){\line(-1, -1){2}}

\put(3.2, 2.5){\footnotesize$u$}
\put(5, 1){\footnotesize$K_{l_{1}}$}
\put(6.8, 1){\footnotesize$K_{l_{2}}$}

\put(10.8, 2.5){\footnotesize$w$}
\put(10.8, 0.5){\footnotesize$y$}
\put(10.8, 4.5){\footnotesize$x$}

\put(8.8, 1.5){\footnotesize$y'$}
\put(8.8, 4.2){\footnotesize$x'$}
\put(13, 3){\footnotesize$z$}

\end{picture}
\caption{The graph $G_{5}(l_{1}, l_{2})$}
\label{G5}
\end{center}
\end{figure}
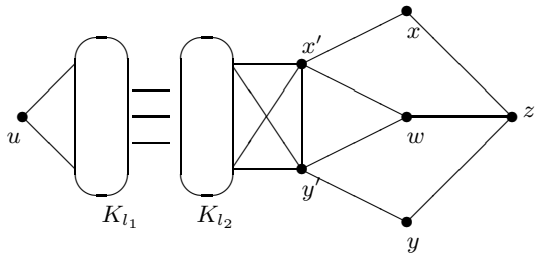
\vskip 10 pt

\indent It is not difficult to show that $G_{5}(l_{1}, l_{2})$ is $5$-$\gamma_{c}$-critical graph. Moreover, $G_{5}(l_{1}, l_{2})$ has $S = \{x', y', z\}$ as a cut set such that $\omega_{o}(G_{5}(l_{1}, l_{2}) - S) = 4 = |S| + 1 > |S| - 1$. By Theorem \ref{thm favaron}, $G_{5}(l_{1}, l_{2})$ is not factor critical.
\vskip 5 pt

\noindent Now, we are ready to prove Theorem \ref{thm 0}. For completeness, we restate the theorem.
\vskip 5 pt

\noindent \textbf{Theorem \ref{thm 0}}.
\emph{Every $k$-$\gamma_{c}$-critical graph of odd order with $\delta \geq 2$ is factor critical if and only if $k \in [2]$. Namely, $\mathcal{Q}(k, 1) = \emptyset$ if and only if $k \in [2]$.}

\proof
Observation \ref{obs 5} implies that if $k = 1$ or $2$, then $\mathcal{Q}(k, 1) = \emptyset$.
\vskip 5 pt

\indent Conversely, Figure~\ref{f:fig001} in Section \ref{pre} yields that $\mathcal{Q}(3, 1) \neq \emptyset$. Lemma \ref{lem 6} yields that $\mathcal{Q}(4, 1) \neq \emptyset$. Moreover, as $G_{5}(l_{1}, l_{2}) \in \mathcal{Q}(5)$, we must have $\mathcal{Q}(5, 1) \neq \emptyset$. We assume that $k \geq 6$. For an odd integer $s \geq 3$, we let $G \in \mathcal{X}(s)$ with $K_{s} = H$ satisfies $(i)$ and $(ii)$ of graphs in the class $\cP(k)$. We moreover let integers $n_{1}, n_{2}, ..., n_{l - 1} \geq 2$ and $n_{l} = 1$ be such that $n_{1} + n_{2} + ... + n_{l - 1} + 1$ is odd number.
\vskip 5 pt

\indent $G(n_{1}, n_{2}, ..., n_{l - 1}, 1) = x_{0} \vee K_{n_{1}} \vee K_{n_{2}} \vee ... \vee K_{n_{l - 1}} \vee K_{n_{l}} \vee \leftidx{_H}{G}$.
\vskip 5 pt

\noindent In view of Theorem \ref{thm gl}, $G(n_{1}, n_{2}, ..., n_{l - 1}, 1)$ is a $(4 + l)$-$\gamma_{c}$-critical graph. Let $V(K_{n_{l}}) = \{y\}$. Thus $G(n_{1}, n_{2}, ..., n_{l - 1}, 1)$ has $S = \{y, a_{1}, a_{2}, ..., a_{s}\}$ as a cut set such that
 \begin{center}
 $\omega_{o}(G(n_{1}, n_{2}, ..., n_{l - 1}, 1) - S) = |S| + 2 > |S| - 1$.
\end{center}
Theorem \ref{thm favaron} then gives that $G(n_{1}, n_{2}, ..., n_{l - 1}, 1)$ is non-factor critical. Therefore $\mathcal{Q}(k, 1) \neq \emptyset$ for all $k \geq 6$. These imply that if $\mathcal{Q}(k, 1) = \emptyset$, then $k = 1$ or $2$. This completes the proof.
\qed

\vskip 5 pt

\indent In view of Theorem \ref{thm 0}, it is natural to think of the bi-criticality of $k$-$\gamma_{c}$-critical graphs with $\delta \geq 3$. Although, we know that if a graph $G$ is not factor critical, then there exists a cut set $S$ such that
\begin{align}
    \omega_{o}(G - S) > |S| - 1 > |S| - 2.\notag
\end{align}
\noindent By Theorem \ref{thm favaron}, regardless with the parity of the orders of graphs, it is most likely there exist $k$-$\gamma_{c}$-critical graphs that are not bi-critical. However, we notice that the graphs that are obtained from the construction in Figure \ref{G5} and Theorem \ref{thm 0} contain a claw, $K_{1, 3}$, as an induced subgraphs. Hence, we may ask if every $k$-$\gamma_{c}$-critical $K_{1, 3}$-free graph with $\delta \geq 3$ is bi-critical. By Theorem \ref{thm chen}, for $1 \leq k \leq 2$, every $k$-$\gamma_{c}$-critical graph is $K_{1, 3}$-free with $\delta \geq 3$ of even order is bi-critical. Therefore, we obtain the following corollary by Observation \ref{obs bicri}
\vskip 5 pt

\begin{Corollary}
If $k \in [2]$, then every $k$-$\gamma_{c}$-critical $K_{1, 3}$-free graph with $\delta \geq 3$ is bi-critical.
\end{Corollary}

\noindent When $k \geq 3$, it turns out that there exist $k$-$\gamma_{c}$-critical graph with $\delta \geq 3$ which is not bi-critical even they do not contain $K_{1, 3}$ as an induced subgraph.
\vskip 5 pt

\noindent Let
\vskip 5 pt

\indent $\tilde{\mathcal{Q}}(k, \ell) :$ the class of $k$-$\gamma_{c}$-critical $K_{1, 3}$-free graphs $G$  with $\delta \geq \ell + 1$ such that\\
\indent $\quad \quad \quad \quad$ $G$ is not $\ell$-factor critical.
\vskip 5 pt

\noindent \textbf{The class $\mathcal{A}(t_{1}, t_{2})$}\\
\indent For an odd number $t_{1} \geq 2$ and an even number $t_{2} \geq 2$, the graph $G$ in this class is obtained from vertices $x_{1}, x_{2}, x_{3}$ and copies of complete graphs $K_{t_{1}}, K_{t_{2}}$ by adding edges according to the join operations:
\begin{itemize}
  \item $x_{1} \vee K_{t_{l}} \vee x_{2} \vee x_{3}$ and
  \item $x_{1} \vee K_{t_{2}} \vee x_{3}$.
\end{itemize}

\begin{Lemma}\label{lem bicri}
If $G \in \mathcal{A}(t_{1}, t_{2})$, then $G \in \mathcal{P}(3) \cap \tilde{\mathcal{Q}}(3, 2)$.
\end{Lemma}
\proof
It is easy to see that all graphs in the class $\mathcal{A}(t_{1}, t_{2})$ is $3$-$\gamma_{c}$-critical and $K_{1, 3}$-free with $\delta \geq 3$. Further, by removing $x_{1}$ and $x_{2}$, the resulting graph has two odd components $K_{t_{1}}$ and $G[V(K_{t_{2}}) \cup \{x_{3}\}]$. Thus, $G$ is not bi-critical.
Therefore, $G \in \tilde{\mathcal{Q}}(3, 2)$. It is easy to see that every graph $G$ in the class $\mathcal{A}(t_{1}, t_{2})$ satisfies the Properties (\emph{i}) and (\emph{ii}) of the class $\mathcal{P}(k)$ by selecting $G[V(K_{t_{2}}) \cup \{x_{3}\}]$ as a maximal complete subgraph. Thus, $G \in \mathcal{P}(3) \cap \tilde{\mathcal{Q}}(3, 2)$ and this completes the proof.
\qed
\vskip 5 pt

\indent Finally,  we will establish the existence of $k$-$\gamma_{c}$-critical $K_{1, 3}$-free graphs with $\delta \geq 3$ which are not  bi-critical for $k \geq 4$ by proving Theorem \ref{thm bicri}. For completeness, we restate theorem.
\vskip 5 pt

\noindent \textbf{Theorem \ref{thm bicri}}. \emph{Every $k$-$\gamma_{c}$-critical $K_{1, 3}$-free graph of even order with $\delta \geq 3$ is bi-critical if and only if $k \in [2]$. Namely, $\tilde{\mathcal{Q}}(k, 2) = \emptyset$ if and only if $k \in [2]$.}

\proof
By Observation \ref{obs bicri}, we have that if $k \in [2]$, then $\tilde{\mathcal{Q}}(k, 2) = \emptyset$.
\vskip 5 pt

\indent Conversely, we will show that if $k \geq 3$, then $\tilde{\mathcal{Q}}(k, 2) \neq \emptyset$. By Lemma \ref{lem bicri}, we have that $\tilde{\mathcal{Q}}(3, 2) \neq \emptyset$. For $k \geq 4$, we let $G \in \mathcal{A}(t_{1}, t_{2})$ such that $t_{1} \geq 3$ is an odd number and $t_{2} \geq 2$ is an even number. Further, we let $x_{0}$ be a vertex and, for $n_{1}, ..., n_{k - 3} \geq 2$ such that $n_{1} + \cdots + n_{k - 3}$ is an odd number, we let $K_{n_{1}}, ..., K_{n - k}$ be $k - 3$ copies of complete graphs of order $n_{1}, ..., n_{k - 3}$, respectively. Let $H = G[V(K_{t_{2}}) \cup \{x_{3}\}]$ and
\vskip 5 pt

\indent $G' = G'(n_{1}, ..., k_{k - 3}, t_{1}, t_{2}) = x_{0} \vee K_{n_{1}} \vee \cdots \vee K_{n_{k - 3}} \vee \leftidx{_H}{G}$.
\vskip 5 pt

\noindent By Theorem \ref{thm gl}, we have that $G'$ is $k$-$\gamma_{c}$-critical graphs. Clearly, $G'$ is $K_{1, 3}$-free and $\delta(G') \geq 3$. By removing $x_{1}$ and $x_{2}$ from $G'$, we have that $G' - x_{1} - x_{2}$ has $K_{t_{1}}$ and $x_{0} \vee K_{n_{1}} \vee \cdots \vee K_{n_{k - 3}} \vee H$ as the two odd components. Therefore,
\vskip 5 pt

\indent $\omega_{o}(G' - \{x_{1}, x_{2}\}) = 2 > 0 = |\{x_{1}, x_{2}\}| - 2$
\vskip 5 pt

\noindent implying that $G'$ is not bi-critical. This completes the proof.
\qed

\end{document}